\documentclass[12pt]{amsart}

\usepackage{etex}
\usepackage{amsmath}
\usepackage{amsfonts}
\usepackage{amssymb}
\usepackage[mathscr]{eucal}
\usepackage[all,cmtip]{xy}
\usepackage{multicol}          %xypic macro for latex2.09
\usepackage{bm}
\usepackage{bbm}
\usepackage{bbding}
\usepackage{graphicx}
\usepackage{txfonts}
\usepackage{amscd}
\usepackage[all]{xy}
\usepackage[shortlabels]{enumitem}
\usepackage{ifpdf}
\ifpdf
  \usepackage[colorlinks,final,backref=page,hyperindex]{hyperref}
\else
  \usepackage[colorlinks,final,backref=page,hyperindex,hypertex]{hyperref}
\fi
\usepackage{tikz}
\usetikzlibrary{positioning}
\usepackage{caption}

\usepackage[T1]{fontenc}
\usepackage{amsmath,amssymb,mathtools}
\usepackage{tikz-cd}
\hypersetup{colorlinks=true,linkcolor=blue,citecolor=blue,urlcolor=blue,
}

\makeatletter

\newtheorem{thm}{Theorem}[section]
\newtheorem{lem}[thm]{Lemma}
\newtheorem{cor}[thm]{Corollary}
\newtheorem{pro}[thm]{Proposition}
\theoremstyle{definition}
\newtheorem{ex}[thm]{Example}
\newtheorem{rmk}[thm]{Remark}
\newtheorem{defi}[thm]{Definition}

\newcommand {\emptycomment}[1]{}

\newcommand{\delete}[1]{{}}

\newcommand{\be }{\begin{equation}}
\newcommand{\ee }{\end{equation}}

\newcommand{\g}{\mathfrak g}

\newcommand{\huaV}{\mathcal{V}}

\newcommand{\huaQ}{\mathcal{Q}}
\newcommand{\huaP}{\mathcal{P}}
\newcommand{\huaT}{\mathcal{T}}

\newcommand{\Id}{{\rm{Id}}}

\newcommand{\br}[1]{   [ \cdot,    \cdot  ]   }

\newcommand{\Hom}{\mathrm{Hom}}
\newcommand{\End}{\mathrm{End}}

\newcommand{\Der}{\mathrm{Der}}

\newcommand{\ABSu}{\mathrm{ABSu}}
\newcommand{\ATSu}{\mathrm{ATSu}}
\newcommand{\ADup}{\mathrm{ADup}}
\newcommand{\ATrip}{\mathrm{ATrip}}

\newcommand{\sgn}{\mathrm{sgn}}

\newcommand{\K}{\mathbf{k}}

\newcommand{\nc}{\newcommand}

\nc{\liead}{{\mathit{LieAd}}}
\nc{\ass}{{\mathit{Ass}}}
\nc{\antiden}{{\mathit{AntiDend}}}
\nc{\antitriden}{{\mathit{AntiTriDend}}}
\nc{\antiprelie}{{\mathit{AntiPreLie}}}
\nc{\antipostlie}{{\mathit{AntiPostLie}}}
\nc{\antipostcom}{{\mathit{AntiPostCom}}}
\nc{\antizinbiel}{{\mathit{AntiZinbiel}}}
\nc{\antiperm}{{\mathit{AntiDupliCom}}}
\nc{\antileibniz}{{\mathit{AntiDupliLie}}}
\nc{\com}{{\mathit{Com}}}
\nc{\antiprepoisson}{{\mathit{AntiPrePoisson}}}
\nc{\poisson}{{\mathit{Poisson}}}
\nc{\comtri}{{\mathit{ComTri}}}
\nc{\triLeib}{{\mathit{TriLeib}}}
\nc{\dend}{{\mathit{Dend}}}
\nc{\dias}{{\mathit{Dias}}}
\nc{\leib}{{\mathit{Leib}}}
\nc{\ldend}{{\mathit{LDend}}}
\nc{\lie}{{\mathit{Lie}}}
\nc{\lquad}{{\mathit{LQuad}}}
\nc{\octo}{{\mathit{Octo}}}
\nc{\perm}{{\mathit{Perm}}}
\nc{\postlie}{{\mathit{PostLie}}}
\nc{\prelie}{{\mathit{PreLie}}}
\nc{\postprelie}{{\mathit{PostPreLie}}}
\nc{\quado}{{\mathit{Quad}}}
\nc{\tridend}{{\mathit{TriDend}}}
\nc{\zinb}{{\mathit{Zinb}}}
\nc{\permtri}{{\mathit{PermTri}}}

\nc{\calt}{{\mathscr T}} \nc{\caltr}{{\mathcal T}^{\,r}}
\nc{\vertset}{\Omega}
\nc{\vin}{{\mathrm{Vin}}}    %decoration set of indices
\nc{\lin}{{\mathrm{Lin}}}    %decoration set of leaves
\nc{\inv}{{\mathrm{In}}}
\nc{\st}{{\mathrm{st}}}
\nc{\ST}{{\mathrm{ST}}}
\nc{\weak}{{\mathrm{weak}}}
\nc{\comp}{{\mathrm{comp}}}
\nc{\COMP}{{\mathrm{COMP}}}
\nc{\WK}{{\mathrm{WK}}}
\nc{\supp}{{\mathrm{supp}}}

\nc{\spann}{{\mathrm{span}}}
\nc{\AntiRB}{{\mathrm{AntiRB}}}

\nc{\gop}{{\,\omega\,}}     % generic binary operation
\nc{\gopb}{{\,\nu\,}}
\nc{\gensp}{V}
\nc{\genbas}{\mathcal{V}}
\nc{\genbi}{\widetilde{V}}
\nc{\gentri}{\widehat{V}}
\nc{\gendup}{\overline{V}}
\nc{\gentrip}{\overrightarrow{V}}

\nc{\svec}[2]{{\textrm{\tiny{$\left(\begin{matrix}#1\\
#2\end{matrix}\right)$}}}}  % column vector
\nc{\cala}{{\mathcal A}}
\nc{\rmi}{{\mathrm{I}}}
\nc{\rmii}{{\mathrm{II}}}
\nc{\rmiii}{{\mathrm{III}}}

\nc{\opd}{\mathcal{P}} %operad P
\nc{\oqd}{\mathcal{Q}}
\nc{\relsp}{R}
\nc{\bvp}{V_P}
\nc{\leafsp}{\mathcal{X}}

\nc{\lir}[1]{\textcolor{red}{Li: #1}}
\nc{\cm}[1]{\textcolor{red}{\underline{CM:}#1 }}
\nc{\wy}[1]{{\marginpar{*}\scriptsize\textcolor{purple}{wy: #1}}}

\begin{document}

\title[Anti-successors and anti-replicators]{Anti-successors and anti-replicators, Manin products and Koszul duality}

\author{Chengming Bai}
\address{Chern Institute of Mathematics and LPMC, Nankai University,
Tianjin 300071, China}
\email{baicm@nankai.edu.cn}

\author{Li Guo}
\address{Department of Mathematics and Computer Science,
         Rutgers University,
         Newark, NJ 07102}
\email{liguo@rutgers.edu}

\author{Yunhe Sheng}
\address{Department of Mathematics, Jilin University, Changchun 130012, Jilin, China}
\email{shengyh@jlu.edu.cn}

\author{You Wang}
\address{Department of Mathematics, Jilin University, Changchun 130012, Jilin, China}
\email{wangyou888@jlu.edu.cn}
%\date{\today}

\begin{abstract}
We develop an operadic theory of anti-splitting and
anti-replicating for binary operads. For a binary operad
$\opd$, we define its anti-bisuccessor $\ABSu(\opd)$
and anti-trisuccessor $\ATSu(\opd)$. Each generating
operation is recovered as the negative sum of two or
three split operations, yielding canonical operad
morphisms $\opd\to\ABSu(\opd)$ and
$\opd\to\ATSu(\opd)$. For binary quadratic operads,
these constructions are described by Manin black
products:
\[
\ABSu(\opd)\cong\antiprelie\bullet\opd,
\qquad
\ATSu(\opd)\cong\antipostlie\bullet\opd.
\]
Dually, we also introduce the canonical quadratic anti-duplicator
and anti-triplicator, together with their combinatorial
counterparts defined by labeled planar binary trees. We show that they coincide precisely under a full-support condition. As applications, we obtain anti-dupli-com and anti-dupli-Lie algebras. We also display an example in which the two constructions differ.
Furthermore, we establish the Koszul duality relations
\[
\widehat{\ADup}(\opd^!)
\cong\ABSu(\opd)^!,
\qquad
\widehat{\ATrip}(\opd^!)
\cong\ATSu(\opd)^!,
\]
providing an effective description of the
canonical quadratic anti-duplicator in terms of
the Manin white product. The framework recovers known constructions, including
anti-dendriform, anti-pre-Lie and anti-pre-Poisson
algebras, and produces further examples and Koszul dual
pairs.
\end{abstract}



\maketitle

\vspace{-1cm}

\tableofcontents

\vspace{-1cm}

\allowdisplaybreaks

\newpage

\section{Introduction}
Splitting and replicating algebraic operations provide two complementary ways to construct
and organize nonassociative structures.  In the ordinary theory these processes are tied to
Manin black and white products and quadratic duality.
Motivated by the recent appearance of anti-pre-Lie and anti-dendriform algebras, we develop
an operadic anti-version of this picture, introducing anti-successors and anti-replicators and
relating them to Manin black and white products.  A central point is that the anti-theory is not obtained by merely changing signs: additional compatibility
relations, and their dual intersections, are essential---and it is precisely these extra relations that force the anti-successor and anti-replicator constructions to be Koszul dual to each other.

\subsection{Splitting and replicating operations}
\label{ss:split}
A basic example of splitting is provided by dendriform algebras: two binary operations are
combined to recover an associative product, and related constructions lead to tridendriform, pre-Lie, pre-Poisson and
other structures; see, among many sources,
\cite{Aguiar2000,BGN,Fard,Loday1,LR}.

The theory of operads provides a natural framework for understanding such constructions in a unified manner~\cite{LV}, with important progress in the study of Koszul duality and Gr\"obner basis methods for operads \cite{DK,GK95}.
It is in the context of operads that the authors of~\cite{BBGN} introduced
the bisuccessor and trisuccessor as splittings of a binary operad and showed,
for binary quadratic operads, that these constructions are
governed by Manin black products and Rota-Baxter operators.  The splitting procedure was later extended from binary operads to general algebraic operads through the configuration
formalism~\cite{PBG2017}. The theory of splitting of operations
has been extended to alternative and Malcev structures
\cite{Madariaga}, and further explored in connection with derived operations and Manin products of operads \cite{KSO}. Subsequently,
pre-Novikov algebras and derived Zinbiel varieties were
investigated in \cite{KMS}, while the behavior of the Dong property under various operadic constructions, including
bisuccessor/trisuccessor (called dendriform splitting) and Manin products, was studied in \cite{KS}.

The complementary process of {\bf replicating} operations has its
motivation from examples such as Leibniz algebras, perm algebras,
dialgebras and trialgebras~\cite{Loday1,LR}. A general theory was
carried out in terms of varieties in universal algebra in
\cite{GubarevKolesnikov2013,KolesnikovVoronin2013}. The terms
\emph{duplicator} and \emph{triplicator} were subsequently used
in~\cite{PBGN} to give a closely related treewise operadic
treatment. Their formulation makes especially transparent three
additional features: the quadratic duality with successors, the
description by Manin white products, and the realization by
averaging operators. A concise overview putting bisuccessors and
duplicators, Manin products, Rota-Baxter operators and averaging
operators into one picture was given earlier in the proceedings
survey \cite{PBGNSurvey2013}. These works are the direct ordinary
counterparts of the constructions developed below.

\subsection{Anti-structures}
In recent years, several new algebraic structures arising from the so-called ``anti-splitting'' of operations have attracted increasing attention. One of the earliest examples is the notion of anti-pre-Lie algebras, arising from  nondegenerate commutative $2$-cocycles on Lie algebras~\cite{LB}. Anti-pre-Lie algebras are also closely related to transposed Poisson algebras~\cite{LB2}. Since then, considerable progress has been made in the study of anti-pre-Lie algebras and related structures~\cite{BGao,ChenZhao}, including compatible anti-pre-Lie algebras~\cite{Zafar} and anti-dendriform algebras~\cite{AANS,BasB}. See \cite{CL} for the geometric structures associated with anti-pre-Lie algebras. Beyond these algebraic constructions, the corresponding anti-Rota-Baxter operators and bialgebra theories have also been developed, which led to analogues of the classical Yang-Baxter equation~\cite{LB,SZ}.

These developments motivate us to construct a general framework,
parallel to the successor and replicator processes described above, that encompasses these new structures in a unified setting.
This is the purpose of the present paper. Note that a different approach to study anti-algebras was proposed in \cite{Ivan} recently. 

\subsection{Main results and organization}

Let \(\mathcal P\) be a binary quadratic operad over a field of characteristic \(0\).  We
introduce its anti-bisuccessor \(\operatorname{ABSu}(\mathcal P)\) and anti-trisuccessor
\(\operatorname{ATSu}(\mathcal P)\) (Definition~\ref{def:ATBSu-sym-operad}).  These constructions recover the known anti-dendriform, anti-pre-Lie and anti-pre-Poisson operads, and lead to anti-post-Lie algebras (Proposition~\ref{pro:anti-post-lie-trisucc-20261104}).

The term \emph{anti-splitting} is justified by reconstruction of the original operation as the \emph{negative} sum of its split components, formalized in the canonical operad morphisms
$\opd\to\ABSu(\opd)$ and
$\opd\to\ATSu(\opd)$
(Theorem~\ref{pro:operad-sub-adj}).  Furthermore, for binary quadratic operads, the two constructions are described by Manin black products
(Theorems~\ref{thm:Manin-product} and~\ref{thm:Manin-product-post}):
\begin{equation}\label{eq:intro-black-v3}
 \operatorname{ABSu}(\mathcal P)
   \cong \operatorname{AntiPreLie}\,\bullet\,\mathcal P,
 \qquad
 \operatorname{ATSu}(\mathcal P)
   \cong \operatorname{AntiPostLie}\,\bullet\,\mathcal P.
\end{equation}

We introduce two versions of the anti-structure of the replicating operations which are related, but serve different purposes. One version is the \emph{canonical} quadratic anti-duplicator $\widehat{\ADup}(\opd)$ and
anti-triplicator $\widehat{\ATrip}(\opd)$, collectively referred to as the canonical quadratic
anti-replicators (Definition~\ref{def-can-rep}). The other consists of the \emph{combinatorial}
counterparts $\ADup(\opd)$ and $\ATrip(\opd)$,
defined via the anti-duplicator and anti-triplicator constructions on labeled planar binary trees (Definition~\ref{def:ADup-sym-operad}).
Thus the combinatorial constructions are more
intuitive and convenient for explicit computations. We show that the combinatorial and canonical quadratic
anti-replicators coincide if and only if the
quadratic relation space has full support with
respect to the chosen basis. More precisely,
their difference is determined by the additional
quadratic relations corresponding to monomials
outside the support of the original relation space (Theorem~\ref{lem:comparison-two-rep}).

Applying the combinatorial construction to  the operads $\com$ and $\lie$, we introduce anti-dupli-com and anti-dupli-Lie algebras (Propositions~\ref{pro:antiperm} and~\ref{pro:antileibniz}). Since the quadratic relation
spaces of $\com$ and $\lie$ have full support, these are also the canonical constructions. On the other hand, for the operad $\mathrm{2Ass}$ of $2$-associative algebras~\cite{Lo08,LR2}, the canonical anti-duplicator and the combinatorial one are different (Example~\ref{ex:2Ass-canonical}).

The advantage of the canonical anti-replicator is that it is in Koszul dual to the anti-splitting construction (Theorems~\ref{thm:koszul-dual} and~\ref{thm:koszul-dual-trip}):
\begin{equation}\label{eq:intro-dual-v3}
\widehat{\ADup}(\opd^!)
\cong\ABSu(\opd)^!,
\qquad
\widehat{\ATrip}(\opd^!)
\cong\ATSu(\opd)^!.
\end{equation}
In particular, we obtain the Koszul dual pairs (Corollaries~\ref{cor:perm-prelie})
\[
\antiperm\cong\antiprelie^!,
\qquad
\antileibniz\cong\antizinbiel^!.
\]
Consequently, the canonical anti-duplicator admits the Manin
white-product description (Theorem~\ref{thm:ADup-antiperm}):
\begin{equation}\label{eq:intro-white-XR-20261104}
\widehat{\ADup}(\opd)\cong\antiperm\bigcirc\opd.
\end{equation}
Here $(-)^!$ denotes the sign-twisted quadratic dual.

The paper is organized as follows. Section~\ref{sec:anti-BSu}
develops anti-bisuccessors and anti-trisuccessors of labeled
planar binary trees and binary operads, together with principal
examples. Section~\ref{sec:manin} describes these constructions in
terms of Manin black products. Section~\ref{sec:anti-rep}
introduces and compares the canonical quadratic and combinatorial
anti-replicators, with illustrations by examples. Finally, Section~\ref{sec:dual} establishes the
Koszul duality between anti-successors and canonical quadratic
anti-replicators and derives the corresponding Manin white-product description.

The processes of splitting and replicating operads
can also be described by Rota-Baxter operators
and averaging operators \cite{BBGN,PBGN}.
For the anti-splitting process, previous studies
of anti-pre-Lie and anti-dendriform algebras
\cite{GLB,LB} have shown their close relationship
with anti-Rota-Baxter operators. A general study
of these operators in connection with the
constructions introduced here, extending these special instances, will be considered separately.

\smallskip

\noindent
{\bf Notation. }
Unless otherwise stated, the
ground field $\K$ is assumed to have characteristic $0$.

\section{The anti-successors of a binary operad}\label{sec:anti-BSu}

We first introduce the notions of the
anti-bisuccessors and anti-trisuccessors of a labeled planar
binary tree.
Then these are applied to define the anti-bisuccessor
and anti-trisuccessor for a nonsymmetric operad and a (symmetric)
operad. Subsequently, we use the computation to recover well-konwn operads and obtain several new operads.

\subsection{The anti-successors of labeled planar binary trees}

\subsubsection{Labeled planar binary trees}

\begin{defi}
\begin{enumerate}
\item[{\rm(i)}]
Let $\calt$ denote the set of planar binary reduced rooted trees together with the trivial tree $\vcenter{\xymatrix@M=4pt@R=8pt@C=4pt{\ar@{-}[d]\\ \\}}$. If $t\in \calt$ has $n$ leaves, we call $t$ an {\bf $n$-tree}. The trivial tree $\vcenter{\xymatrix@M=4pt@R=8pt@C=4pt{\ar@{-}[d]\\ \\}}$ has one leaf.
\item[{\rm(ii)}]
 Let $\vertset$ be a set. By a {\bf decorated tree} we mean a tree $t$ of $\calt$ together with a decoration on the vertices of $t$ by elements of $\vertset$ and a decoration on the leaves of $t$ by distinct positive integers. Let $t(\vertset)$ denote the set of decorated trees for $t$ and denote by
$
\calt(\vertset):=\coprod_{ t \in \calt} t(\vertset).
$
If $\tau\in t(\vertset)$ for an $n$-tree $t$, we call $\tau$ a {\bf labeled $n$-tree}.
\item[{\rm(iii)}]
For $\tau\in \calt(\vertset)$, we let $\vin(\tau)$ (resp. $\lin(\tau)$) denote the set (resp. ordered set) of labels of the vertices (resp. leaves) of $\tau$.
\item[{\rm(iv)}]
Let $\tau_{l},\tau_{r}\in \calt(\vertset)$ with disjoint sets of leaf labels and $\omega\in \vertset$. The {\bf grafting of $\tau_{l}$ and $\tau_{r}$ along} $\omega$ is denoted by $\tau=
  \tau_{l} \vee_{\omega} \tau_{r}$, which is an element in $\calt(\vertset)$.
\item[{\rm(v)}]
 For $\tau\in \calt(\vertset)$ with $|\lin(\tau)|>1$, we let $\tau= \tau_{l} \vee_{\omega} \tau_{r}$ denote the unique decomposition of $\tau$ as a grafting of $\tau_{l}$ and $\tau_{r}$
  in $\calt(\vertset)$ along $\omega \in \vertset$.
\end{enumerate}
\end{defi}

Let $t$ be a tree and let $V$ be a graded vector space.  concentrated in arity $2$, so that $V=V(2)$. Recall that the free nonsymmetric operad $\mathcal{T}_{\hspace*{-0.1cm}ns}(\gensp)$ on $V$ is given by the vector space
$
\mathcal{T}_{\hspace*{-0.1cm}ns}(\gensp) := \displaystyle{ \oplus_{t \in \calt } } \ t[\gensp],
$
where $t[\gensp]$ is the treewise tensor module associated to $t$. This module is explicitly given by
$
t[\gensp] := \displaystyle{\otimes_{v \in \vin(t)}} \gensp_{|\inv(v)|},
$
where $\inv(v)$ denotes the incoming edges of $v$, labeled by $1,\cdots,|\inv(v)|$ from left to right. See \cite{LV} for more details. A basis $\genbas$ of $\gensp$ induces a basis $t(\genbas)$ of $t[\gensp]$ and a basis $\calt(\genbas)$ of $\mathcal{T}_{\hspace*{-0.1cm}ns}(\gensp)$. Consequently, any element of $t[\gensp]$ can be represented as a linear combination of elements in $t(\genbas)$.

\subsubsection{Anti-bisuccessors of labeled planar binary trees}

\begin{defi}\label{def:anti-bisu}{\rm
Let $\gensp$ be a vector space with a basis $\genbas$.
\begin{enumerate}
\item[{\rm(i)}]
Define a vector space
\begin{equation}\label{genbi}
\genbi= \gensp \otimes (\mathbf{k} \prec \oplus~ \mathbf{k} \succ),
\end{equation}
where we denote $(\omega~\otimes \prec)$ (resp. $(\omega~\otimes \succ)$) by
$\begin{pmatrix}
  \omega \\
  \prec
  \end{pmatrix}$
(resp. $\begin{pmatrix}
  \omega \\
  \succ
  \end{pmatrix}$)
for $\omega\in \gensp$. Then $\genbas \times \{ \prec,\succ\}$ is a basis of $\genbi$. $\genbi$ can be viewed as a graded module concentrated in arity $2$.
\item[{\rm(ii)}]
For a labeled $n$-tree $\tau$ in $\calt(\genbas)$, define $\widetilde{\tau}$ in $\mathcal{T}_{\hspace*{-0.1cm}ns}(\genbi)$ as follows:
\begin{itemize}
 \item[$\bullet$] $  \widetilde{|}=|$
 \item[$\bullet$] when $n\geq 2$, $\widetilde{\tau}$ is obtained by replacing each decoration $\gop \in \vin(\tau)$ by
     $$\begin{pmatrix}
  \gop \\
  \ast
  \end{pmatrix}:=-\begin{pmatrix}
  \gop \\
  \prec
  \end{pmatrix}-\begin{pmatrix}
  \gop \\
  \succ
  \end{pmatrix}.$$
\end{itemize}
We extend this definition to $\mathcal{T}_{\hspace*{-0.1cm}ns}(\genbi)$ by linearity.
\end{enumerate}
}\end{defi}

\begin{defi}
Let $\gensp$ be a vector space with a basis $\genbas$. Let $\tau$ be a labeled $n$-tree in
$\calt(\genbas)$. The {\bf anti-bisuccessor} $\ABSu_{x}(\tau)$ of $\tau$ with respect to a leaf $x\in \lin(\tau)$ is an element of $\mathcal{T}_{\hspace*{-0.1cm}ns}(\genbi)$ defined by induction on
$|\lin(\tau)|$ as follows:
\begin{itemize}
 \item[$\bullet$] $\ABSu_{x}(~|~)=|$~;
 \item[$\bullet$] Let $k\geq 1$. Assume that $\ABSu_{x}(\tau)$ has been defined for $\tau$ with $|\lin(\tau)|\leq k$.  Then for a labeled $(k+1)$-tree $\tau\in \calt(\genbas)$ with its decomposition $\tau_{l} \vee_{\gop} \tau_{r}$, we define
{\small     $$
     \ABSu_{x}(\tau)= \ABSu_{x}(\tau_{l} \vee_{\gop} \tau_{r}) = \left\{\begin{array}{ll}  \ABSu_{x}(\tau_{l}) \vee_{\begin{pmatrix}
  \gop \\
  \prec
  \end{pmatrix}} \widetilde{\tau_{r}},\quad & \text{if}~x \in \lin(\tau_{l});  \\
     \widetilde{\tau_{l}} \vee_{\begin{pmatrix}
  \gop \\
  \succ
  \end{pmatrix}} \ABSu_{x}(\tau_{r}), \quad & \text{if}~x \in \lin(\tau_{r}). \end{array}\right.
     $$
}
\end{itemize}
\end{defi}

The following proposition characterizes the anti-bisuccessors of labeled trees explicitly.

\begin{pro}\label{pro:ABSu}
Let $\gensp$ be a vector space with a basis $\genbas$, $\tau$ be a
labeled $n$-tree in $\calt(\genbas)$ and $x$ be in $\lin(\tau)$.
The anti-bisuccessor $\ABSu_{x}(\tau)$ of $\tau$ is obtained by
relabelling each vertex of $\tau$ according to the following
rules:
\begin{enumerate}
\item[{\rm(a)}]
We replace the label $\gop$ of each vertex on the path from the root to the leaf $x$ of $\tau$ by
\begin{itemize}
 \item[{\rm(i)}] $\begin{pmatrix}
  \gop \\
  \prec
  \end{pmatrix}$,\quad if the path turns left at this vertex;
\item[\rm{(ii)}] $\begin{pmatrix}
  \gop \\
  \succ
  \end{pmatrix}$,\quad if the path turns right at this vertex;
\end{itemize}
\item[{\rm(b)}]
we replace the label $\gop$ of each vertex not on the path from the root to the leaf $x$ of $\tau$ by
$$\begin{pmatrix}
  \gop \\
  \ast
  \end{pmatrix}:=-\begin{pmatrix}
  \gop \\
  \prec
  \end{pmatrix}-\begin{pmatrix}
  \gop \\
  \succ
  \end{pmatrix}.$$
\end{enumerate}
\end{pro}

\begin{proof}
This follows by induction on $|\lin(\tau)|\geq 1$. We omit the details for brevity.
\end{proof}

\begin{ex}
{\footnotesize$$\ABSu_{x_{2}} \left( \begin{xy}
(0,8)*+{x_{1}}="v1";(10,8)*+{x_{2}}="v2";(20,8)*+{x_{3}}="v3"; (30,8)*+{x_{4}}="v4";
(5,0)*+{\omega_{1}}="v5";(25,0)*+{\omega_{2}}="v6"; (15,-6)*+{\omega_{3}}="v7";(15,-12)*+{}="v8";
{\ar@{-}"v5";"v1"};{\ar@{-}"v5";"v2"};{\ar@{-}"v6";"v3"};{\ar@{-}"v6";"v4"};{\ar@{-}"v7";"v5"};{\ar@{-}"v7";"v6"};{\ar@{-}"v8";"v7"}
\end{xy} \right)= \begin{xy}
(0,10)*+{x_{1}}="v1";(10,10)*+{x_{2}}="v2";(20,10)*+{x_{3}}="v3"; (30,10)*+{x_{4}}="v4";
(5,0)*+{\svec{\omega_{1}}{\succ}}="v5";(25,0)* +{\svec{\omega_{2}}{\ast}}="v6"; (15,-6)*+{\svec{\omega_{3}}{\prec}}="v7";(15,-15)*+{}="v8";
{\ar@{-}"v5";"v1"};{\ar@{->}"v5";"v2"};{\ar@{-}"v6";"v3"};{\ar@{-}"v6";"v4"};{\ar@{->}"v7";"v5"};{\ar@{-}"v7";"v6"};{\ar@{->}"v8";"v7"}
\end{xy}
=-\begin{xy}
(0,10)*+{x_{1}}="v1";(10,10)*+{x_{2}}="v2";(20,10)*+{x_{3}}="v3"; (30,10)*+{x_{4}}="v4";
(5,0)*+{\svec{\omega_{1}}{\succ}}="v5";(25,0)* +{\svec{\omega_{2}}{\prec}}="v6"; (15,-6)*+{\svec{\omega_{3}}{\prec}}="v7";(15,-15)*+{}="v8";
{\ar@{-}"v5";"v1"};{\ar@{-}"v5";"v2"};{\ar@{-}"v6";"v3"};{\ar@{-}"v6";"v4"};{\ar@{-}"v7";"v5"};{\ar@{-}"v7";"v6"};{\ar@{-}"v8";"v7"}
\end{xy}-\begin{xy}
(0,10)*+{x_{1}}="v1";(10,10)*+{x_{2}}="v2";(20,10)*+{x_{3}}="v3"; (30,10)*+{x_{4}}="v4";
(5,0)*+{\svec{\omega_{1}}{\succ}}="v5";(25,0)* +{\svec{\omega_{2}}{\succ}}="v6"; (15,-6)*+{\svec{\omega_{3}}{\prec}}="v7";(15,-15)*+{}="v8";
{\ar@{-}"v5";"v1"};{\ar@{-}"v5";"v2"};{\ar@{-}"v6";"v3"};{\ar@{-}"v6";"v4"};{\ar@{-}"v7";"v5"};{\ar@{-}"v7";"v6"};{\ar@{-}"v8";"v7"}
\end{xy}
$$}
\end{ex}

\begin{lem}\label{group-act}
Let $\gensp$ be a vector space with a basis $\genbas$, $\tau$ be a labeled $n$-tree in
$\calt(\genbas)$ and $x$ be in $\lin(\tau)$. Then the following relation holds
$$ \ABSu_{\sigma^{-1}(x)}(\tau^{\sigma})=\ABSu_{x}(\tau) ^{\sigma},\quad \forall \sigma \in \mathbb{S}_n.   $$
\end{lem}

\begin{proof}
This follows by inspecting the symmetric-group action on a tree.
\end{proof}

\subsubsection{Anti-trisuccessors of labeled planar binary trees}

\begin{defi}\label{def:anti-trisu}
Let $\gensp$ be a vector space with basis $\genbas$.
\begin{enumerate}
\item[{\rm(i)}]
Define a vector space
\begin{eqnarray}\label{gentri}
\gentri= \gensp \otimes (\mathbf{k} \prec \oplus~ \mathbf{k}\succ\oplus~\mathbf{k}~\cdot~),
\end{eqnarray}
where we denote $(\omega~\otimes \prec)$ (resp. $(\omega~\otimes \succ)$, resp. $(\omega~\otimes\cdot)$) by
$\begin{pmatrix}
  \omega \\
  \prec
  \end{pmatrix}$
(resp. $\begin{pmatrix}
  \omega \\
  \succ
  \end{pmatrix}$, resp. $\begin{pmatrix}
  \omega \\
  \cdot
  \end{pmatrix}$
  )
for $\omega\in \gensp$. Then $\genbas \times \{ \prec,\succ,\cdot \}$ is a basis of $\gentri$. $\gentri$ can be viewed as an arity-graded vector space concentrated in arity $2$.
\item[{\rm(ii)}]
For a labeled $n$-tree $\tau$ in $\calt(\genbas)$, define $\widehat{\tau}$ in $\mathcal{T}_{\hspace*{-0.1cm}ns}(\gentri)$ as follows:
\begin{itemize}
 \item[$\bullet$] $  \widehat{|}=|$
 \item[$\bullet$] when $n\geq 2$, $\widehat{\tau}$ is obtained by replacing each decoration $\gop \in \vin(\tau)$ by
     $$\begin{pmatrix}
  \gop \\
  \star
  \end{pmatrix}:=-\begin{pmatrix}
  \gop \\
  \prec
  \end{pmatrix}-\begin{pmatrix}
  \gop \\
  \succ
  \end{pmatrix}-\begin{pmatrix}
  \gop \\
  \cdot
  \end{pmatrix}.$$
\end{itemize}
We extend this definition to $\mathcal{T}_{\hspace*{-0.1cm}ns}(\gentri)$ by linearity.
\end{enumerate}
\end{defi}

\begin{defi}
Let $\gensp$ be a vector space with a basis $\genbas$. Let $\tau$ be a labeled $n$-tree in
$\calt(\genbas)$ and $J$ be a nonempty subset of $\lin(\tau)$. The {\bf anti-trisuccessor} $\ATSu_{J}(\tau)$ of $\tau$ with respect to $J$ is an element of $\mathcal{T}_{\hspace*{-0.1cm}ns}(\gentri)$ defined by induction on
$|\lin(\tau)|$ as follows:
\begin{itemize}
 \item[$\bullet$] $\ATSu_{J}(~|~)=|$~;
 \item[$\bullet$] Assume that $\ATSu_{J}(\tau)$ has been defined for $\tau$ with $|\lin(\tau)|\leq k$. Then for a labeled $(k+1)$-tree $\tau\in \calt(\genbas)$ with its decomposition $\tau_{l} \vee_{\gop} \tau_{r}$, we define
     $$
     \ATSu_{J}(\tau)= \ATSu_{J}(\tau_{l} \vee_{\gop} \tau_{r}) = \left\{\begin{array}{lll}  \ATSu_{J}(\tau_{l}) \vee_{\begin{pmatrix}
  \gop \\
  \prec
  \end{pmatrix}} \widehat{\tau_{r}},\quad & \text{if}~J \subseteq \lin(\tau_{l}),  \\
     \widehat{\tau_{l}} \vee_{\begin{pmatrix}
  \gop \\
  \succ
  \end{pmatrix}} \ATSu_{J}(\tau_{r}), \quad & \text{if}~J \subseteq \lin(\tau_{r}), \\
  \ATSu_{J\cap \lin(\tau_{l})}(\tau_{l}) \vee_{\begin{pmatrix}
  \gop \\
  \cdot
  \end{pmatrix}}
  \ATSu_{J\cap \lin(\tau_{r})}(\tau_{r}), \quad & {\mathrm{otherwise.}}
   \end{array}\right.
     $$
\end{itemize}
\end{defi}

The following proposition gives an explicit characterization of the anti-trisuccessors.

\begin{pro}
Let $\gensp$ be a vector space with a basis $\genbas$, $\tau$ be a
labeled $n$-tree in $\calt(\genbas)$ and $J$ be a nonempty subset
of $\lin(\tau)$. The anti-trisuccessor $\ATSu_{J}(\tau)$ is
obtained by relabelling each vertex of $\tau$ according to the
following rules:
\begin{enumerate}
\item[{\rm(a)}]
We replace the label $\gop$ of each vertex on at least one of the paths from the root to the leaves in $J$ by
\begin{itemize}
 \item[\rm{(i)}] $\begin{pmatrix}
  \gop \\
  \prec
  \end{pmatrix}$, if all such paths turn left at this vertex;
\item[\rm{(ii)}] $\begin{pmatrix}
  \gop \\
  \succ
  \end{pmatrix}$, if all such paths turn right at this vertex;
\item[\rm{(iii)}] $\begin{pmatrix}
  \gop \\
  \cdot
  \end{pmatrix}$, if some of these paths turn left and some of these paths turn right at this vertex;
\end{itemize}
\item[{\rm(b)}]
We replace the label $\gop$ of each other vertex by
$$\begin{pmatrix}
  \gop \\
  \star
  \end{pmatrix}:=-\begin{pmatrix}
  \gop \\
  \prec
  \end{pmatrix}-\begin{pmatrix}
  \gop \\
  \succ
  \end{pmatrix}-\begin{pmatrix}
  \gop \\
  \cdot
  \end{pmatrix}.$$
\end{enumerate}
\end{pro}

\begin{ex}
{\footnotesize$$\ATSu_{\{x_{1},x_{4}\}} \left( \begin{xy}
(0,8)*+{x_{1}}="v1";(10,8)*+{x_{2}}="v2";(20,8)*+{x_{3}}="v3"; (30,8)*+{x_{4}}="v4";
(5,0)*+{\omega_{1}}="v5";(25,0)*+{\omega_{2}}="v6"; (15,-6)*+{\omega_{3}}="v7";(15,-12)*+{}="v8";
{\ar@{-}"v5";"v1"};{\ar@{-}"v5";"v2"};{\ar@{-}"v6";"v3"};{\ar@{-}"v6";"v4"};{\ar@{-}"v7";"v5"};{\ar@{-}"v7";"v6"};{\ar@{-}"v8";"v7"}
\end{xy} \right)= \begin{xy}
(0,10)*+{x_{1}}="v1";(10,10)*+{x_{2}}="v2";(20,10)*+{x_{3}}="v3"; (30,10)*+{x_{4}}="v4";
(5,0)*+{\svec{\omega_{1}}{\prec}}="v5";(25,0)* +{\svec{\omega_{2}}{\succ}}="v6"; (15,-6)*+{\svec{\omega_{3}}{\cdot}}="v7";(15,-15)*+{}="v8";
{\ar@{->}"v5";"v1"};{\ar@{-}"v5";"v2"};{\ar@{-}"v6";"v3"};{\ar@{->}"v6";"v4"};{\ar@{->}"v7";"v5"};{\ar@{->}"v7";"v6"};{\ar@{->}"v8";"v7"}
\end{xy}
$$}
\end{ex}

We have the following property of anti-trisuccessors under permutations.

\begin{lem}\label{group-act-2}
Let $\gensp$ be a vector space with a basis $\genbas$, $\tau$ be a labeled $n$-tree in
$\calt(\genbas)$ and $J$ be a nonempty subset of $\lin(\tau)$. Then the following relation holds
$$ \ATSu_{\sigma^{-1}(J)}(\tau^{\sigma})=\ATSu_{J}(\tau) ^{\sigma},\quad \forall \sigma \in \mathbb{S}_n.   $$
\end{lem}

\begin{proof}
By inspection of the action of the symmetric group on a tree.
\end{proof}

\subsection{The anti-successors of a binary operad}

We first recall some preliminaries and properties of a binary quadratic operad for later application. See \cite{LV,Val2} for more details.

An {\bf $\mathbb{S}$-module} over $\mathbf{k}$ is a family
  $ M:=\{ M(n)\}_{n \geq 0}=\{M(0),M(1),\dots,M(n),\dots\} $
  of right $\mathbf{k} \mathbb{S}_{n}$-modules $M(n)$, where $\mathbb{S}_{n}$ is the symmetric group for all $n \geq 0$. If every $M(n)$ is the trivial $\mathbf{k} \mathbb{S}_{n}$-module, then $M$ is called an {\bf $\mathbb{N}$-module}.
\begin{defi}
 
 {\rm(a)} A nonsymmetric operad is an $\mathbb{N}$-module $\opd=\{ \opd(n)\}_{n \geq 0}$ equipped with {\bf partial compositions}:
    $$ \circ_{i}:=\circ_{m,n,i}: \opd(m) \otimes \opd(n) \to \opd(m+n-1),\quad 1\leq i\leq m$$
    such that the following equalities hold for all $\eta\in \opd(l),\mu \in \opd(m)$ and $\nu\in \opd(n)$:
    \begin{enumerate}
  \item[{\rm(i)}] $(\eta \circ_i \mu) \circ_{i-1+j} \nu= \eta \circ_i (\mu \circ_j \nu), \quad 1\leq i\leq l,1\leq j\leq m.$
  \item[{\rm(ii)}]  $(\eta \circ_i \mu) \circ_{k-1+m} \nu= (\eta \circ_k \nu)\circ_i \mu,\quad 1 \leq i <k \leq l.$
  \item[{\rm(iii)}] There is an element $\Id \in \opd(1)$ such that $\mu \circ \Id=\mu$ and $\Id \circ \mu=\mu$.
     \end{enumerate}
     
 {\rm(b)} A nonsymmetric operad $\opd=\{ \opd(n)\}_{n \geq 0}$ is called {\bf a (symmetric) operad}, if $\opd=\{ \opd(n)\}_{n \geq 0}$ is an $\mathbb{S}$-module and satisfies the following additional conditions
    \begin{enumerate}
  \item[{\rm(iv)}] For any $\sigma \in \mathbb{S}_{m}$, we have
  $$ \mu^{\sigma} \circ_i \nu=(\mu \circ_{\sigma(i)} \nu)^{\sigma'},\quad 1 \leq i\leq  m, $$
  where $\sigma'\in \mathbb{S}_{m+n-1}$ restricts to the bijection from $\{1,\dots,m+n-1\} \backslash \{i,\dots,i+n-1\}$ to $\{1,\dots,m+n-1\} \backslash \{\sigma(i),\dots,\sigma(i)+n-1\}$ according to
  $\sigma$ and restricts to the bijection from $\{i,\dots,i+m-1\}$ to $\{\sigma(i),\dots,\sigma(i)+m-1\}$ according to the identity (in the natural order).
  \item[{\rm(v)}] For any $\tau \in \mathbb{S}_{n}$, we have
  $ \mu \circ_i \nu^{\tau}=(\mu \circ_{i} \nu)^{\tau'},\quad 1 \leq i\leq  m, $
  where $\tau'\in \mathbb{S}_{m+n-1}$ is the permutation which acts by the identity, except on the block $\{i,\dots,i+n-1\}$ on which it acts via $\tau$.
  \end{enumerate}
 
\end{defi}

A {\bf morphism of operads} from $\opd$ to $\oqd$ is a family of $\mathbb{S}_{n}$-equivariant maps $f_n:\opd(n) \to \oqd(n)$ that is compatible with the compositions of operads. An operad
$\opd$ is called {\bf reduced} if $\opd(0)=0$. Throughout this paper, all operads are assumed to be reduced.

A standard example of an operad is the {\bf endomorphism operad}  $\End_{V}$ for a vector space $V$, given by
$ \End_{V}(n):=\Hom(V^{\otimes n},V). $ Let $\opd$ be an operad. A {\bf $\opd$-algebra} is a vector space $V$ with a morphism of operads $\rho:\opd \to \End_{V}$. In this case, we say that the operad $\opd$ {\bf encodes} $\opd$-algebras.

A class of operads can be given more explicitly as follows \cite[Section 7.6.2]{LV}.
Let $V=\{V(n)\}_{n \geq 0}$ be an $\mathbb{S}$-module with $V(1)=\mathbf{k} \Id$. If $V(k)=0$ for all $k\neq 1,2$, we call $V$ {\bf binary}. In this case, we denote $V:=V(2)$. Thus an arbitrary binary operad is a quotient
$$ \opd(V,R):= \huaT(V)/(R)$$
for some $\mathbb{S}_2$-module $V$ and an operadic ideal $(R)$ of $\huaT(V)$ generated by a subset $R\subseteq \huaT(V)$. $V$ is called the {\bf space of generators} and $R$ is called the {\bf relation space}. If all elements of $R$ are quadratic, namely linear combinations of the form $\omega \circ_{i} \nu$ for $\omega,\nu \in V, i=1,2$, then the operad  $ \opd(V,R)$ is called {\bf binary quadratic}.

To determine the possible quadratic relations for $V$, we note that
$$ \huaT(V)(3)=(V \otimes_{\mathbb{S}_2} (V \otimes \mathbf{k} \oplus \mathbf{k} \otimes V))\otimes_{\mathbb{S}_2} \mathbf{k}[\mathbb{S}_3],    $$
which can be identified with $3$ copies of $V \otimes V$, denoted by $V \circ_{\rmi} V, V \circ_{\rmii} V$ and $V \circ_{\rmiii} V$, following the convention in \cite{Val2}. Then as a vector space, $\huaT(V)(3)$ is generated by elements of the form
\begin{equation}\label{eq:i-ii-iii}
\omega \circ_{\rmi} \nu (\leftrightarrow  (x \nu y)\omega z),\quad \omega \circ_{\rmii} \nu (\leftrightarrow  (y \nu z)\omega x),\quad \omega \circ_{\rmiii} \nu (\leftrightarrow  (z \nu x)\omega y), \quad \forall \omega,\nu \in V.
\end{equation}

For an operad whose space of binary generators is
$V(2):=\mathbf{k}\mu\oplus \mathbf{k}\mu'$, regarded as a right ${\mathbb S}_2$-module via
$\mu\cdot(12)=\mu'$ (equivalently, $\mu'\cdot(12)=\mu$), we will adopt the convention in \cite[p. 129]{Val2} and denote the twelve elements of $ \huaT(V)(3)$ by $v_i~(1\leq i \leq 12)$ in the following table.
\begin{equation}
\label{eq:table}
\begin{tabular}{|c|l|c|l|c|l|}
\hline
\(v_1\) & \(\mu \circ_{\rmi} \mu \leftrightarrow (xy)z\) & \(v_5\) & \(\mu \circ_{\rmiii} \mu \leftrightarrow (zx)y\) & \(v_9\) & \(\mu \circ_{\rmii} \mu \leftrightarrow (yz)x\) \\
\hline
\(v_2\) & \(\mu' \circ_{\rmii} \mu \leftrightarrow x(yz)\) & \(v_6\) & \(\mu' \circ_{\rmi} \mu \leftrightarrow z(xy)\) & \(v_{10}\) & \(\mu' \circ_{\rmiii} \mu \leftrightarrow y(zx)\) \\
\hline
\(v_3\) & \(\mu' \circ_{\rmii} \mu' \leftrightarrow x(zy)\) & \(v_7\) & \(\mu' \circ_{\rmi} \mu' \leftrightarrow z(yx)\) & \(v_{11}\) & \(\mu' \circ_{\rmiii} \mu' \leftrightarrow y(xz)\) \\
\hline
\(v_4\) & \(\mu \circ_{\rmiii} \mu' \leftrightarrow (xz)y\) & \(v_8\) & \(\mu \circ_{\rmii} \mu' \leftrightarrow (zy)x\) & \(v_{12}\) & \(\mu \circ_{\rmi} \mu' \leftrightarrow (yx)z\) \\
\hline
\end{tabular}
\end{equation}

Let $\gensp=\gensp(2)$ be an $\mathbb{S}$-module concentrated in arity $2$ with a linear basis $\genbas$. For any finite set $\leafsp$ of cardinal $n$, we can define the following coinvariant space
$$  \gensp(\leafsp):=\begin{pmatrix} \bigoplus\limits_{f:\b{n}\to \leafsp} \gensp(n)  \end{pmatrix}_{\mathbb{S}_n}, $$
where the sum is over all the bijections from $\b{n}:=\{1,\cdots,n\}$ to $\leafsp$ and the symmetric group $\mathbb{S}_n$ acts diagonally.

Denote by $\mathbb{T}$ the set of isomorphism classes of reduced binary trees \cite[Appendix C]{LV}. For $\mathrm{t}\in \mathbb{T}$, define $\mathrm{t}[\gensp]$ as the treewise tensor $\mathbb{S}$-module associated to $\mathrm{t}$, which is explicitly given by
$
\mathrm{t}[\gensp] := \displaystyle{\otimes_{v \in \vin(\mathrm{t})}} \gensp(\inv(v)).
$
For further details, see \cite[Section 5.5.1]{LV}. Then the free operad $\mathcal{T}(\gensp)$ on an $\mathbb{S}$-module $\gensp=\gensp(2)$ is given by the $\mathbb{S}$-module
$
\mathcal{T}(\gensp) := \displaystyle{ \oplus_{\mathrm{t} \in \mathbb{T} } } \ \mathrm{t}[\gensp].
$
Each tree $\mathrm{t}\in \mathbb{T}$ can be represented by a planar tree $t$ in $\calt$ by choosing a total order on the set of inputs of each vertex of $\mathrm{t}$. We can define an equivalence relation $\thicksim$ on $\calt$ as follows: two planar binary trees in $\calt$ are equivalent if they represent the same tree in $\mathbb{T}$. It gives rise to a bijection $\mathbb{T}\cong \calt/\thicksim$. Furthermore, for any planar binary tree $t$ in $\calt$ which represents $\mathrm{t}\in \mathbb{T}$, we have $t[\gensp]\cong\mathrm{t}[\gensp]$ according to \cite[Section 2.8]{Hoffbeck}. Thus we have
$
\mathcal{T}(\gensp) \cong \displaystyle{ \oplus_{t \in \mathfrak{R} } } \ t[\gensp],
$
where $\mathfrak{R}$ is a set of representatives of $\calt/\thicksim$ and $\mathfrak{R}\cong \mathbb{T}$. It is straightforward to check that
$\mathfrak{R}(\genbas):=\{ \tau\in t(\genbas)| t \in \mathfrak{R}\}$ is a linear basis of the free operad $\mathcal{T}(\gensp)$.

We now introduce the anti-successors of a binary symmetric operad.

\begin{defi}\label{def:ATBSu-sym-operad}
Let $\opd=\mathcal{T}(\gensp)/ (\relsp)$ be a binary (symmetric) operad on the $\mathbb{S}$-module $\gensp=\gensp(2)$ concentrated in arity $2$ with a $\mathbf{k}[\mathbb{S}_2]$-basis $\genbas$, such that the space of relations $\relsp$ is spanned as an $\mathbb{S}$-module, by locally homogeneous elements of the form
\begin{equation}\label{eq:relation-1}
R:=\Biggl\{  r_s:=\sum_{i} c_{s,i} \tau_{s,i}~|~ c_{s,i} \in \mathbf{k}, \tau_{s,i} \in  \{t(\genbas),t\in \mathfrak{R}\}, 1 \leq s\leq k, k\geq 1 \Biggr\},
\end{equation}
where $\mathfrak{R}$ is a set of representatives of $\calt/\thicksim$ defined as above. For each fixed $s$, all trees $\tau_{s,i}$ have the same set of leaves $\lin(\tau_{s,i})$.

\begin{enumerate}
\item[{\rm(i)}]
The {\bf anti-bisuccessor} of $\opd$ is defined to be the binary operad
$$\ABSu(\opd):=\mathcal{T}(\genbi)/ (\ABSu(\relsp),\st(\relsp)),$$
where the $\mathbb{S}_2$-action on $\genbi$ is given by
$$ \begin{pmatrix}
  \gop \\
  \prec
  \end{pmatrix}^{(12)}:=\begin{pmatrix}
  \gop^{(12)} \\
  \succ
  \end{pmatrix},\quad
  \begin{pmatrix}
  \gop \\
  \succ
  \end{pmatrix}^{(12)}:=\begin{pmatrix}
  \gop^{(12)} \\
  \prec
  \end{pmatrix}, \quad
  \forall \gop\in \gensp.
$$
The space of relations is generated as an $\mathbb{S}$-module, given by
\begin{eqnarray*} \ABSu(\relsp)&:=&\Biggl\{ \ABSu_x(r_s)=\sum_{i} c_{s,i} \ABSu_x(\tau_{s,i})~|~ x\in \lin(\tau_{s,i}), 1 \leq s \leq k \Biggr\},\\
 \st(\relsp)&:=&\Biggl\{  \sum_{i} c_{s,i} \st(\tau_{s,i})=\sum_{i} c_{s,i} \widetilde{\tau}_{s,i}~|~1 \leq s \leq k   \Biggr\}.  
\end{eqnarray*}
Here $\st(\tau_{s,i})=\widetilde{\tau}_{s,i}$ is obtained by replacing each vertex decoration $\gop$ of $\tau_{s,i}$ by
$$\begin{pmatrix}
  \gop \\
  \ast
  \end{pmatrix}:=-\begin{pmatrix}
  \gop \\
  \prec
  \end{pmatrix}-\begin{pmatrix}
  \gop \\
  \succ
  \end{pmatrix}.$$
\item[{\rm(ii)}]
The {\bf anti-trisuccessor} of $\opd$ is defined to be the binary operad
$$\ATSu(\opd):=\mathcal{T}(\gentri)/ (\ATSu(\relsp),\ST(\relsp)),$$
where the $\mathbb{S}_2$-action on $\gentri$ is given by
$$ \begin{pmatrix}
  \gop \\
  \prec
  \end{pmatrix}^{(12)}:=\begin{pmatrix}
  \gop^{(12)} \\
  \succ
  \end{pmatrix},\quad
  \begin{pmatrix}
  \gop \\
  \succ
  \end{pmatrix}^{(12)}:=\begin{pmatrix}
  \gop^{(12)} \\
  \prec
  \end{pmatrix}, \quad
    \begin{pmatrix}
  \gop \\
  \cdot
  \end{pmatrix}^{(12)}:=\begin{pmatrix}
  \gop^{(12)} \\
  \cdot
  \end{pmatrix},\quad
  \forall \gop\in \gensp.
$$
The space of relations is generated as an $\mathbb{S}$-module, given by
\begin{eqnarray*}\ATSu(\relsp)&:=&\Biggl\{ \ATSu_{J}(r_s)=\sum_{i} c_{s,i} \ATSu_{J}(\tau_{s,i})~\Biggr| \varnothing \neq J \subseteq \lin(\tau_{s,i}), 1 \leq s \leq k \Biggr\},\\
   \ST(\relsp)&:=&\Biggl\{  \sum_{i} c_{s,i} \ST(\tau_{s,i})=\sum_{i} c_{s,i} \widehat{\tau}_{s,i}~\Biggr|~1 \leq s \leq k   \Biggr\}. 
   \end{eqnarray*}
Here $\ST(\tau_{s,i})=\widehat{\tau}_{s,i}$ is obtained by replacing each vertex decoration $\gop$ of $\tau_{s,i}$ by
$$\begin{pmatrix}
  \gop \\
  \star
  \end{pmatrix}:=-\begin{pmatrix}
  \gop \\
  \prec
  \end{pmatrix}-\begin{pmatrix}
  \gop \\
  \succ
  \end{pmatrix}-\begin{pmatrix}
  \gop \\
  \cdot
  \end{pmatrix}.$$
\end{enumerate}
\end{defi}

An argument analogous to the proof of Proposition~2.20 in \cite{BBGN} shows that the anti-successor of a binary operad $\opd=\mathcal{T}(\gensp)/ (\relsp)$ depends neither on the linear basis $\genbas$ of $\gensp$ nor on the set $\mathfrak{R}$.

\begin{rmk}
The symbols $\star$ and $\ast$ record the signed total operations underlying the anti-successor construction. In contrast with an ordinary splitting, the original operation will be recovered from the split operations through the \emph{negative} total sum. Theorem~\ref{pro:operad-sub-adj} below makes this precise by showing that these signed total operations satisfy the defining relations of the original operad.
\end{rmk}

\begin{rmk}
The above definitions also apply to binary nonsymmetric operads by omitting the symmetric group actions and replacing the free symmetric operad $\mathcal{T}(\gensp)$ with the free nonsymmetric operad $\mathcal{T}_{\hspace*{-0.1cm}ns}(\gensp)$. More precisely, let
$
\opd=\mathcal{T}_{\hspace*{-0.1cm}ns}(\gensp)/(\relsp),
$
where $\gensp=\gensp(2)$ is a vector space with a linear basis $\genbas$, and $\relsp$ is spanned by locally homogeneous elements of the form
$$
r_s=\sum_i c_{s,i}\tau_{s,i},\qquad
c_{s,i}\in\mathbf{k},\quad \tau_{s,i}\in\calt(\genbas),
\quad 1\leq s\leq k.
$$
For each fixed $s$, all trees $\tau_{s,i}$ have the same arity $n_s$. The anti-bisuccessor and anti-trisuccessor are then defined by the same formulas:
$$
\ABSu(\opd):=
\mathcal{T}_{\hspace*{-0.1cm}ns}(\genbi)/
(\ABSu(\relsp),\st(\relsp)),
\qquad
\ATSu(\opd):=
\mathcal{T}_{\hspace*{-0.1cm}ns}(\gentri)/
(\ATSu(\relsp),\ST(\relsp)),
$$
where the relations are understood as nonsymmetric operadic relations. In this setting, both constructions are independent of the choice of a basis of the vector space of generating operations, by linearity of the anti-successor constructions and the treewise tensor structure of the free nonsymmetric operad.
\end{rmk}

We now justify the terminology \emph{anti-splitting}. Although the
defining operations are obtained by splitting each generator into
two or three pieces, the operation that reconstructs the original
$\opd$-structure is the negative total operation introduced in
Definition~\ref{def:ATBSu-sym-operad}. The following theorem
expresses this fact intrinsically at the operadic level.

\begin{thm}\label{pro:operad-sub-adj}
Let $\opd=\mathcal{T}(\gensp)/ (\relsp)$ be a binary operad.
\begin{enumerate}
\item  \label{i:osa1}
There is a morphism of operads from $\opd$ to $\ABSu(\opd)$ which extends the linear map from $\gensp$ to $\widetilde{V}$ defined by
\begin{equation}\label{opd-ABSuopd}
\omega \to \begin{pmatrix}
  \omega\\
  \ast
  \end{pmatrix},\quad \forall \omega\in \gensp.
\end{equation}
\item \label{i:osa2}
There is a morphism of operads from $\opd$ to $\ATSu(\opd)$ which extends the linear map from $\gensp$ to $\widehat{V}$ defined by
\begin{equation}\label{opd-ATSuopd-1}
\omega \to \begin{pmatrix}
  \omega\\
  \star
  \end{pmatrix},\quad \forall \omega\in \gensp.
\end{equation}
\end{enumerate}
\end{thm}

\begin{proof}
Assume that $R$ is given by Eq. \eqref{eq:relation-1}.

\eqref{i:osa1} It is easy to check that the linear map defined in Eq.
\eqref{opd-ABSuopd} is $\mathbb{S}_2$-equivariant and hence induces a morphism of
operads from $\mathcal{T}(\gensp)$ to $\ABSu(\opd)$. Since $\sum_i
c_{s,i}\widetilde{\tau}_{s,i}\in \st(R)$, there exists a morphism of operads from $\opd$ to $\ABSu(\opd)$.

\eqref{i:osa2}  Since we have $\sum_i c_{s,i}\widehat{\tau}_{s,i}\in \ST(R)$, the proof is similar to the proof of Item \eqref{i:osa1}.
\end{proof}

\begin{rmk}\label{rmk:sub-adj}
Note that a morphism of operads $f:\opd\to\oqd$ induces a functor from the category of $\oqd$-algebras to the category of $\opd$-algebras. More precisely, if $(V,\rho:\oqd\to \End_V)$ is a $\oqd$-algebra, then the composite
$ \opd \xrightarrow{f} \oqd \xrightarrow{\rho} \End_V $
makes $V$ into a $\opd$-algebra.
\end{rmk}
\begin{rmk}
Theorem~\ref{pro:operad-sub-adj} gives an intrinsic explanation of
the term \emph{anti-splitting}. Indeed, the morphism
$\opd\to\ABSu(\opd)$ sends each generator $\omega$ to
$\omega_\star=-(\omega_\prec+\omega_\succ)$. Hence every
$\ABSu(\opd)$-algebra canonically becomes a $\opd$-algebra under
the negative total operation. Likewise, every
$\ATSu(\opd)$-algebra carries the original $\opd$-structure
through $\omega_\star=-(\omega_\prec+\omega_\succ+\omega_\cdot)$.
Thus the construction is a genuine splitting procedure, with the
prefix ``anti-'' reflected precisely in the negative
reconstruction of the original operation. Note that, for any binary quadratic operad $\opd$, $(A,\cdot)$ is a $\opd$-algebra if and only if $(A,-\cdot)$ is a $\opd$-algebra. Therefore, for notational simplicity, we replace $\ast$ and $\star$ by $-\ast$ and $-\star$, respectively, in the following examples, while retaining the notation $\ast$ and $\star$ for the resulting operations.
\end{rmk}

\subsection{Examples of anti-successors of binary operads}

\subsubsection{Examples of anti-bisuccessors}

Recall that an {\bf anti-dendriform algebra} \cite{GLB} is a vector space $A$ equipped with two bilinear operations $\lhd,\rhd$ satisfying the following equations:
\begin{eqnarray}
\label{anti-den-5} &&x\rhd (y\rhd z)=-(x \ast y)\rhd z=-x\lhd (y \ast z)=(x\lhd y)\lhd z;\\
\label{anti-den-6} &&(x\rhd y)\lhd z=x\rhd (y \lhd z),
\end{eqnarray}
for all $x,y,z\in A$, where $x \ast y:=x \rhd y+ x \lhd y$.

\begin{pro}\label{pro:antiden}
The operad $\antiden$ encoding anti-dendriform algebras is the anti-bisuccessor of the nonsymmetric operad $\ass$ encoding associative algebras.
\end{pro}

\begin{proof}
Let $\cdot$ denote the product of $\ass$.
 Using the abbreviations $\rhd:=\begin{pmatrix}
  \cdot\\
  \succ
  \end{pmatrix}$, $\lhd:=\begin{pmatrix}
  \cdot\\
  \prec
  \end{pmatrix}$ and $\ast:=\rhd + \lhd$, by a direct calculation, we have
\begin{eqnarray}
\label{anti-den-1} \ABSu_{x}((x \cdot y)\cdot z- x\cdot (y\cdot z)) &=& (x\lhd y)\lhd z+ x\lhd (y \ast z);\\
\label{anti-den-2} \ABSu_{y}((x \cdot y)\cdot z- x\cdot (y\cdot z)) &=& (x\rhd y)\lhd z- x\rhd (y \lhd z);\\
\label{anti-den-3} \ABSu_{z}((x \cdot y)\cdot z- x\cdot (y\cdot z)) &=& -(x \ast y)\rhd z- x\rhd (y \rhd z);\\
\label{anti-den-4} \st((x \cdot y)\cdot z- x\cdot (y\cdot z)) &=& (x \ast y)\ast z-x\ast (y \ast z),
\end{eqnarray}
giving the same relations \eqref{anti-den-5}-\eqref{anti-den-6} of anti-dendriform algebras introduced in \cite{GLB}. Therefore, the operad $\antiden$ is the anti-bisuccessor of the operad $\ass$.
\end{proof}

Recall that a {\bf (left) anti-pre-Lie algebra} \cite{LB} is a vector space $\g$ equipped with a bilinear operation $\circ$ satisfying the following relations
\begin{eqnarray}
\label{anti-prelie-1} x \circ (y \circ z)-y \circ (x \circ z)+[x,y]_{\circ} \circ z &=&0;\\
\label{anti-prelie-2} [[x,y]_{\circ},z]_{\circ}+[[y,z]_{\circ},x]_{\circ}+[[z,x]_{\circ},y]_{\circ}&=&0, \quad x, y, z\in \g,
\end{eqnarray}
where $[x,y]_{\circ}:= x \circ y-y \circ x$.

\begin{pro}\label{pro:anti-pre-lie}
The operad $\antiprelie$ encoding anti-pre-Lie algebras is the anti-bisuccessor of the operad $\lie$ encoding Lie algebras.
\end{pro}

\begin{proof}
Let $\mu$ be the generating operation of the operad $\lie$. Set $\circ:=\begin{pmatrix}
  \mu\\
  \succ
  \end{pmatrix}$. Since $\begin{pmatrix}
  \mu\\
  \prec
  \end{pmatrix}^{(12)}=\begin{pmatrix}
  \mu^{(12)}\\
  \succ
  \end{pmatrix}=-\begin{pmatrix}
  \mu\\
  \succ
  \end{pmatrix}=-~\circ$, we have $\begin{pmatrix}
  \mu\\
  \prec
  \end{pmatrix}=-\circ^{(12)}.$ The space of the relations of $\lie$ is generated as an $\mathbb{S}_3$-module by
$$ v_1+v_5+v_9= \mu \circ_{\rmi} \mu+\mu \circ_{\rmii} \mu+\mu \circ_{\rmiii} \mu. $$
Then we have
\begin{eqnarray*}
\ABSu_{x}(v_1+v_5+v_9)&=& z \circ (y \circ x)-[y,z]_{\circ} \circ x-y \circ (z \circ x);\\
\ABSu_{y}(v_1+v_5+v_9)&=& x \circ (z \circ y)-[z,x]_{\circ} \circ y-z \circ (x \circ y);\\
\ABSu_{z}(v_1+v_5+v_9)&=& y \circ (x \circ z)-[x,y]_{\circ} \circ z-x \circ (y \circ z);\\
\st(v_1+v_5+v_9)&=& [[x,y]_{\circ},z]_{\circ}+[[y,z]_{\circ},x]_{\circ}+[[z,x]_{\circ},y]_{\circ},
\end{eqnarray*}
which coincide with the relations \eqref{anti-prelie-1} and
\eqref{anti-prelie-2}
of anti-pre-Lie algebras. Therefore,
$\ABSu(\lie)=\antiprelie$.
\end{proof}

Now we recall the notions of anti-Zinbiel algebras and
anti-Poisson algebras introduced in \cite{LB2023} in the study of
new splittings of Poisson algebras, and show that the operads
$\antizinbiel$ and $\antiprepoisson$ encoding anti-Zinbiel
algebras and anti-pre-Poisson algebras are anti-bisuccessors of
the operads $\com$ and $\poisson$ encoding commutative
associative algebras and Poisson algebras respectively.

\begin{defi}{\rm(\cite{LB2023})}
An {\bf anti-Zinbiel algebra} $(A,\diamond)$ is a vector space $A$ equipped with a bilinear operation $\diamond:A \otimes A \to A$ satisfying the following equation
\begin{equation}\label{anti-zinbiel-5}
x \diamond (y \diamond z)=-(x \diamond y+ y\diamond x)\diamond z=z \diamond (y \diamond x),\quad \forall x,y,z\in A.
\end{equation}
\end{defi}

\begin{defi}{\rm(\cite{LB2023})}\label{def:anti-pre-poisson}
An {\bf anti-pre-Poisson algebra} $(A,\diamond,\circ)$ is a vector space $A$ equipped with two bilinear operations $\diamond,\circ:A \otimes A \to A$
such that $(A,\diamond)$ is an anti-Zinbiel algebra and $(A,\circ)$ is an anti-pre-Lie algebra and for all $x,y,z\in A,$ the following conditions hold:
\begin{eqnarray}
\label{anti-pre-poisson-1} (y \diamond z+z \diamond y)\circ x+z \diamond (y \circ x)+y \diamond (z \circ x)&=&0;\\
\label{anti-pre-poisson-2} (x \circ z- z\circ x)\diamond y+x \circ (z \diamond y)- z \diamond (x \circ y)&=&0;\\
\label{anti-pre-poisson-3} (y \diamond z+z \diamond y)\circ x+(x \circ z- z\circ x)\diamond y+(x \circ y- y\circ x)\diamond z&=&0.
\end{eqnarray}
\end{defi}

\begin{pro}\label{pro:antizinbiel}
\begin{enumerate}
  \item \label{i:az1} The operad $\antizinbiel$  is the anti-bisuccessor of the operad $\com$;

   \item \label{i:az2} The operad $\antiprepoisson$  is the anti-bisuccessor of the operad $\poisson$.
\end{enumerate}
\end{pro}

\begin{proof}
We only prove \eqref{i:az2}; Item \eqref{i:az1} is similar.
Recall that the operad $\poisson:=\huaT(V)/(R_{Poisson})$ is a binary operad on the $\mathbb{S}$-module $V=V(2)$ concentrated in arity $2$ with a $\mathbf{k}[\mathbb{S}_2]$-basis $\huaV:=\{\omega,\mu\}$, such that the space of relations $R_{Poisson}$ is generated as an $\mathbb{S}_3$-module by
\[
\left\{
\begin{array}{l}
\omega \circ_{\rmi} \omega-\omega \circ_{\rmii} \omega;\\
\mu \circ_{\rmi} \mu+\mu \circ_{\rmii} \mu+\mu \circ_{\rmiii} \mu;\\
\mu' \circ_{\rmii} \omega-\omega \circ_{\rmi} \mu-\omega' \circ_{\rmiii} \mu'.
\end{array}
\right.
\]
By a direct calculation, we have
\begin{eqnarray*}
\ABSu_{x}(\mu' \circ_{\rmii} \omega-\omega \circ_{\rmi} \mu-\omega' \circ_{\rmiii} \mu')&=& (y \diamond z+z \diamond y)\circ x+z \diamond (y \circ x)+y \diamond (z \circ x);\\
\ABSu_{y}(\mu' \circ_{\rmii} \omega-\omega \circ_{\rmi} \mu-\omega' \circ_{\rmiii} \mu')&=& x \circ (z \diamond y)- z \diamond (x \circ y)+(x \circ z- z\circ x)\diamond y;\\
\ABSu_{z}(\mu' \circ_{\rmii} \omega-\omega \circ_{\rmi} \mu-\omega' \circ_{\rmiii} \mu')&=& x \circ (y \diamond z)+(x \circ y- y\circ x)\diamond z- y \diamond (x \circ z);\\
\st(\mu' \circ_{\rmii} \omega-\omega \circ_{\rmi} \mu-\omega' \circ_{\rmiii} \mu') &=& \{x, y\cdot z\}-\{x, y\}\cdot z-y \cdot \{x, z\},
\end{eqnarray*}
where $\{x, y\}:= x \circ y-y \circ x$ and $ x\cdot y:=x \diamond y+y \diamond x$. It is easy to check that these correspond to the relations \eqref{anti-pre-poisson-1}--\eqref{anti-pre-poisson-3} in Definition \ref{def:anti-pre-poisson}. Therefore, we have $\ABSu(\poisson)=\antiprepoisson$.
\end{proof}

\subsubsection{Examples of anti-trisuccessors}
We first introduce the notion of anti-tridendriform algebras and show that the operad $\antitriden$ encoding anti-tridendriform algebras is the anti-trisuccessor of the operad $\ass$.

\begin{defi}
An {\bf anti-tridendriform algebra} $(A,<,>,\cdot)$ is a vector space $A$ equipped with an associative product $\cdot$ and two bilinear operations $<,>:A \otimes A \to A$ satisfying the following equations
\begin{eqnarray}
\label{anti-den-tri-1} &&(x < y)< z=-x < (y \star z), \quad (x > y)< z=x > (y < z),\quad (x > y)> z=-(x \star y)> z,\\
\label{anti-den-tri-2} && (x \cdot y)< z=x \cdot (y < z), \quad (x > y)\cdot z=x > (y \cdot z), \quad (x < y)\cdot z=x \cdot (y > z),\\
\label{anti-den-tri-3} &&(x \star y)\star z=x \star (y \star z),
\end{eqnarray}
for all $x,y,z\in A$, where the binary operation $\star:A \otimes A \to A$ is given by $x \star y=x\cdot y+x > y+x < y$.
\end{defi}

\begin{rmk}
Let $(A,<,>,\cdot)$ be an anti-tridendriform algebra. If the associative algebra $(A,\cdot)$ is trivial, then $(A,<,>)$ becomes an anti-dendriform algebra. Therefore, an anti-tridendriform algebra can be viewed as a generalization of an anti-dendriform algebra.
\end{rmk}

By Theorem \ref{pro:operad-sub-adj} and Remark \ref{rmk:sub-adj},
we obtain the following result.

\begin{pro}\label{pro:sub-adj-anti-triden}
Let $(A,<,>,\cdot)$ be an anti-tridendriform algebra.
%Then we have \cm{delete ``Then we have" or change to be ``Then we have the following results."}
\begin{enumerate}
  \item\label{it:it1} $(A,\star)$ is an associative algebra, which is called the {\bf sub-adjacent associative algebra} of $(A,<,>,\cdot)$ and denoted by $A^c$, where the multiplication
      $\star:A \otimes A \to A$ is given by
      $$x \star y=x\cdot y+x > y+x < y, \quad \forall x,y \in A.$$
  \item\label{it:it2} $(-L_>,-R_<):A \to \End(A)$ gives an action of the associative algebra $A^c$ on $A$, where the left and right multiplication operators $L_>,R_<:A \to \End(A)$ are defined by $L_>(x)(y)=x > y, R_<(x)(y)=y < x$ respectively, for all $x,y\in A$.
\end{enumerate}
\end{pro}

\begin{proof} For item~(\ref{it:it1}), the assertion follows directly from Theorem
\ref{pro:operad-sub-adj} and Remark \ref{rmk:sub-adj}.

For item~(\ref{it:it2}), the assertion follows by a straightforward verification.
\end{proof}

\begin{pro}\label{pro:anti-den-tri}
The operad $\antitriden$ is the anti-trisuccessor of the operad $\ass$.
\end{pro}

\begin{proof}
Let $\nu$ be the generating operation of the operad $\ass$ satisfying the relation $v_1-v_2=\nu \circ_{\rmi} \nu-\nu' \circ_{\rmii} \nu$. Using the abbreviations $>:=\begin{pmatrix}
  \nu\\
  \succ
  \end{pmatrix}$, $<:=\begin{pmatrix}
  \nu\\
  \prec
  \end{pmatrix}$, $\cdot:=\begin{pmatrix}
  \nu\\
  \cdot
  \end{pmatrix}$ and $\star:=>+<+~\cdot$, by a direct calculation, we have
\begin{eqnarray}
\label{eq:anti-den-tri-1} \ATSu_{\{x\}}(v_1-v_2)&=& (x < y)< z+x < (y \star z);\\
\label{eq:anti-den-tri-2} \ATSu_{\{y\}}(v_1-v_2)&=& (x > y)< z-x > (y < z);\\
\label{eq:anti-den-tri-3} \ATSu_{\{z\}}(v_1-v_2)&=& -(x > y)> z-(x \star y)> z;\\
\label{eq:anti-den-tri-4} \ATSu_{\{x,y\}}(v_1-v_2)&=& (x \cdot y)< z-x \cdot (y < z);\\
\label{eq:anti-den-tri-5} \ATSu_{\{x,z\}}(v_1-v_2)&=& (x > y)\cdot z-x > (y \cdot z);\\
\label{eq:anti-den-tri-6} \ATSu_{\{y,z\}}(v_1-v_2)&=& (x < y)\cdot z-x \cdot (y > z);\\
\label{eq:anti-den-tri-7} \ATSu_{\{x,y,z\}}(v_1-v_2)&=& (x \cdot y)\cdot z-x \cdot (y \cdot z);\\
\label{anti-den-tri-ST} \ST(v_1-v_2)&=& (x \star y)\star z-x \star (y \star z),
\end{eqnarray}
which coincide with the relations
of anti-tridendriform algebras. Thus, the operad $\antitriden$ is
the anti-trisuccessor of the operad $\ass$.
\end{proof}

The notion of post-Lie algebras was introduced in \cite{Val1} in the study of operads and has found various applications in geometric numerical analysis and regularity structures. Now we introduce the notion of anti-post-Lie algebras and show that the operad $\antipostlie$ is the anti-trisuccessor of the operad $\lie$.

\begin{defi}
An {\bf anti-post-Lie algebra} $(\g,[\cdot,\cdot]_\g,\triangleright)$ consists of a Lie algebra $(\g,[\cdot,\cdot]_\g)$ and a binary operation $\triangleright:\g \otimes \g \to \g$ such that for all $x,y,z\in \g$,
\begin{eqnarray}
\label{Post-1} &&x\triangleright [y,z]_\g=[x\triangleright y,z]_\g+[y,x\triangleright z]_\g,\\
\label{Post-2} &&([x,y]_\g+x\triangleright y-y\triangleright x) \triangleright z=-x\triangleright(y\triangleright z)+y\triangleright(x\triangleright z),\\
\label{Post-3} &&([x,y]_\g+x\triangleright y-y\triangleright x) \triangleright z+([y,z]_\g+y\triangleright z-z\triangleright y) \triangleright x+([z,x]_\g+z\triangleright x-x\triangleright z) \triangleright y=0.
\end{eqnarray}

\end{defi}

\begin{rmk}
Let $(\g,[\cdot,\cdot]_\g,\triangleright)$ be an anti-post-Lie algebra. If the Lie bracket $(\g,[\cdot,\cdot]_\g)$ vanishes, then $(\g,\triangleright)$ degenerates to an anti-pre-Lie algebra.
Thus, an anti-post-Lie algebra can be viewed as a nonabelian generalization of an anti-pre-Lie algebra.
\end{rmk}

\begin{pro}\label{pro:sub-adj-anti-postlie}
Any anti-post-Lie algebra $(\g,[\cdot,\cdot]_\g,\triangleright)$ gives rise to a new Lie algebra $(\g,[\cdot,\cdot]_\triangleright)$ given by
\begin{eqnarray}\label{eq:sub-adj-Lie}
 [x,y]_\triangleright:=x \triangleright y-y \triangleright x+[x,y]_\g,\quad \forall x,y\in \g,
\end{eqnarray}
which is called the {\bf sub-adjacent Lie algebra}  and denoted by $\g_{\triangleright}$. Moreover, $-L_{\triangleright}:\g \to \Der(\g)$ defined by $-L_{\triangleright}(x)(y)=-x \triangleright y$ is an action of $\g_{\triangleright}$ on the Lie algebra $(\g,[\cdot,\cdot]_\g)$.
\end{pro}

\begin{proof} The first assertion follows directly from Theorem
\ref{pro:operad-sub-adj} and Remark \ref{rmk:sub-adj}. The second assertion follows by a straightforward verification.
\end{proof}

\begin{pro}
\label{pro:anti-post-lie-trisucc-20261104}
The operad $\antipostlie$ encoding anti-post-Lie algebras is the anti-trisuccessor of the operad $\lie$.
\end{pro}

\begin{proof}
Let $\mu$ be the generating operation of the operad $\lie$. Set $\triangleright:=\begin{pmatrix}
  \mu\\
  \succ
  \end{pmatrix}$ and $[\cdot,\cdot]:=\begin{pmatrix}
  \mu\\
  \cdot
  \end{pmatrix}$. Since $\begin{pmatrix}
  \mu\\
  \prec
  \end{pmatrix}^{(12)}=\begin{pmatrix}
  \mu^{(12)}\\
  \succ
  \end{pmatrix}=-\begin{pmatrix}
  \mu\\
  \succ
  \end{pmatrix}=-\triangleright$, we have
   $\begin{pmatrix}
  \mu\\
  \prec
  \end{pmatrix}=-\triangleright^{(12)}.$ Furthermore, we also have
\begin{eqnarray*}
\ATSu_{\{x\}}(v_1+v_5+v_9)&=& z \triangleright (y \triangleright x)-y \triangleright (z \triangleright x)+(z\triangleright y-y \triangleright z-[y,z])\triangleright x;\\
\ATSu_{\{y\}}(v_1+v_5+v_9)&=& x \triangleright (z \triangleright y)-z \triangleright (x \triangleright y)+(x\triangleright z-z \triangleright x-[z,x])\triangleright y;\\
\ATSu_{\{z\}}(v_1+v_5+v_9)&=& y \triangleright (x \triangleright z)-x \triangleright (y \triangleright z)+(y\triangleright x-x \triangleright y-[x,y])\triangleright z;\\
\ATSu_{\{x,y\}}(v_1+v_5+v_9)&=& -z \triangleright [x,y]-[z\triangleright y,x]+[z \triangleright x,y];\\
\ATSu_{\{x,z\}}(v_1+v_5+v_9)&=& -y \triangleright [z,x]-[y\triangleright x,z]+[y \triangleright z,x];\\
\ATSu_{\{y,z\}}(v_1+v_5+v_9)&=& -x \triangleright [y,z]-[x\triangleright z,y]+[x \triangleright y,z];\\
\ATSu_{\{x,y,z\}}(v_1+v_5+v_9)&=& [[x,y],z]+
[[y,z],x]+[[z,x],y];\\
\ST(v_1+v_5+v_9)&=& [[x,y]_\triangleright,z]_\triangleright+
[[y,z]_\triangleright,x]_\triangleright+[[z,x]_\triangleright,y]_\triangleright,
\end{eqnarray*}
which agree with the relations
\eqref{Post-1}-\eqref{Post-3} of anti-post-Lie algebras.
Therefore, $\ATSu(\lie)=\antipostlie$.
\end{proof}

Now we introduce the notion of anti-post-com algebras and show that the operad $\antipostcom$ encoding anti-post-com algebras is the anti-trisuccessor of the operad $\com$.

\begin{defi}
An {\bf anti-post-com algebra} $(A,\cdot,\rhd)$ consists of a commutative associative algebra $(A,\cdot)$ and a binary operation $\rhd:A \otimes A \to A$ such that
\begin{eqnarray}
\label{Post-com-1} -(x \star y) \rhd z &=& x \rhd (y \rhd z);\\
\label{Post-com-2} z \rhd (x \cdot y) &=& x \cdot (z \rhd y);\\
\label{Post-com-3} x \rhd (y \rhd z) &=& z \rhd (y \rhd x),
\end{eqnarray}
where the bilinear product $\star:A \otimes A \to A$ is given by $x \star y:=x \rhd y+y \rhd x+x \cdot y$, for all $x,y,z\in A$.
\end{defi}

\begin{rmk}
Let $(A,\cdot,\rhd)$ be an anti-post-com algebra. If the commutative associative algebra $(A,\cdot)$ is trivial, then $(A,\rhd)$ degenerates to an anti-Zinbiel algebra.
Thus, an anti-post-com algebra can be viewed as a generalization of an anti-Zinbiel algebra.

Moreover, let $(A,<,>,\cdot)$ be an anti-tridendriform algebra. If the associative algebra $(A,\cdot)$ is commutative and the two operations $<,>$ satisfy $x < y=y > x$, then $(A,\cdot,>)$ is an anti-post-com algebra.
\end{rmk}

\begin{pro}\label{pro:sub-adj-anti-postcom}
Any anti-post-com algebra $(A,\cdot,\rhd)$ gives rise to a new commutative associative algebra $(A,\star)$ given by
\begin{eqnarray}\label{eq:sub-adj-com}
x \star y:=x \rhd y+y \rhd x+x \cdot y, \quad \forall x,y\in A,
\end{eqnarray}
which is called the {\bf sub-adjacent commutative (associative) algebra} of the anti-post-com algebra $(A,\cdot,\rhd)$ and denoted by $A_{\rhd}$. Moreover, $-L_{\rhd}:A \to \End(A)$ defined by $-L_{\rhd}(x)(y)=-x \rhd y$ is an action of $A_{\rhd}$ on the commutative associative algebra $(A,\cdot)$.
\end{pro}

\begin{proof} The first assertion follows directly from Theorem
\ref{pro:operad-sub-adj} and Remark \ref{rmk:sub-adj}. The second assertion follows by a straightforward verification.
\end{proof}

Let $\antipostcom$ be the operad encoding anti-post-com algebras.

\begin{pro}\label{pro:anti-postcom}
The operad $\antipostcom$  is the anti-trisuccessor of the operad $\com$.
\end{pro}

\begin{proof}
Let $\omega$ be the generating operation of the operad $\com$. Set $\rhd:=\begin{pmatrix}
  \omega\\
  \succ
  \end{pmatrix}$ and $\cdot:=\begin{pmatrix}
  \omega\\
  \cdot
  \end{pmatrix}$. Since $\begin{pmatrix}
  \omega\\
  \prec
  \end{pmatrix}^{(12)}=\begin{pmatrix}
  \omega^{(12)}\\
  \succ
  \end{pmatrix}=\begin{pmatrix}
  \omega\\
  \succ
  \end{pmatrix}=\rhd$, we have $\begin{pmatrix}
  \omega\\
  \prec
  \end{pmatrix}=\rhd^{(12)}.$ The space of relations of $\com$ is generated as an $\mathbb{S}_3$-module by
$$ v_1-v_9= \omega \circ_{\rmi} \omega-\omega \circ_{\rmii} \omega. $$
Then we have
{\small
\begin{eqnarray}
\label{anti-post-com-1} \ATSu_{\{x\}}(v_1-v_9)&=& z \rhd (y \rhd x)+(y \star z)\rhd x;\\
\label{anti-post-com-2} \ATSu_{\{y\}}(v_1-v_9)&=& z \rhd (x \rhd y)-x \rhd (z \rhd y);\\
\label{anti-post-com-3} \ATSu_{\{z\}}(v_1-v_9)&=& -(x \star y)\rhd z-x \rhd (y \rhd z);\\
\label{anti-post-com-4} \ATSu_{\{x,y\}}(v_1-v_9)&=& z \rhd (x \cdot y)-(z \rhd y) \cdot x;\\
\label{anti-post-com-5} \ATSu_{\{x,z\}}(v_1-v_9)&=& (y \rhd x)\cdot z-(y \rhd z)\cdot x;\\
\label{anti-post-com-6} \ATSu_{\{y,z\}}(v_1-v_9)&=& (x \rhd y)\cdot z- x \rhd (y \cdot z);\\
\label{anti-post-com-7} \ATSu_{\{x,y,z\}}(v_1-v_9)&=& (x \cdot y)\cdot z-(y \cdot z)\cdot x;\\
\label{anti-post-com-st} \ST(v_1-v_9)&=& (x \star y)\star z-(y \star z)\star x,
\end{eqnarray}
}
which corresponds to the relations in Eqs. \eqref{Post-com-1}--\eqref{Post-com-3} of anti-post-com algebras. Thus, the operad $\antipostcom$ is the anti-trisuccessor of the operad $\com$.
\end{proof}

\section{Anti-successors and Manin product}
\label{sec:manin}
We identify the anti-trisuccessor of a binary quadratic operad $\opd$ with the
Manin black product of $\antipostlie$ with $\opd$. Similarly, we
also prove that there is an operadic isomorphism between the
anti-bisuccessor of a binary quadratic operad $\opd$ and the Manin
black product of $\antiprelie$ with $\opd$.

\begin{defi}{\rm(\cite{GK95,Val2})}
Let $\huaP=\huaT(V)/(R)$ and $\huaQ=\huaT(W)/(S)$ be two binary quadratic operads with finite-dimensional generating spaces. Define their {\bf Manin black product} by the formula
$$ \huaP \bullet \huaQ:= \huaT(V \otimes W \otimes \mathbf{k}. \sgn_{\mathbb{S}_2})/(\Psi(R \otimes S)), $$
where $\Psi$ is defined in \cite[Section 4.3]{Val2}.
\end{defi}

\begin{thm}\label{thm:Manin-product-post}
Let $\opd$ be a binary quadratic operad. There is an isomorphism of operads
$$ \ATSu(\opd) \cong \antipostlie \bullet \opd. $$
\end{thm}

\begin{proof}
Denote the generating operations $[\cdot,\cdot]$ and $\triangleright$ of $\antipostlie$ by $\beta$ and $\varepsilon$ respectively. Then we have $\beta'=-\beta$. The space of relations of $\antipostlie$ is generated as a vector space by
\begin{eqnarray*}
 &&\beta \circ_{\rmi} \beta +\beta  \circ_{\rmii} \beta +\beta  \circ_{\rmiii} \beta,\\
 &&-\varepsilon' \circ_{\rmi} \beta+\beta \circ_{\rmiii} \varepsilon+\beta' \circ_{\rmii} \varepsilon',\\
 &&-\varepsilon' \circ_{\rmii} \beta+\beta \circ_{\rmi} \varepsilon+\beta' \circ_{\rmiii} \varepsilon',\\
 &&-\varepsilon' \circ_{\rmiii} \beta+\beta \circ_{\rmii} \varepsilon+\beta' \circ_{\rmi} \varepsilon',\\
 &&-\varepsilon \circ_{\rmii} \beta-\varepsilon \circ_{\rmii} \varepsilon+\varepsilon \circ_{\rmii} \varepsilon'-\varepsilon'\circ_{\rmiii} \varepsilon+\varepsilon'\circ_{\rmi} \varepsilon',\\
 &&-\varepsilon \circ_{\rmiii} \beta-\varepsilon \circ_{\rmiii} \varepsilon+\varepsilon \circ_{\rmiii} \varepsilon'-\varepsilon'\circ_{\rmi} \varepsilon+\varepsilon'\circ_{\rmii} \varepsilon',\\
 &&-\varepsilon \circ_{\rmi} \beta-\varepsilon \circ_{\rmi} \varepsilon+\varepsilon \circ_{\rmi} \varepsilon'-\varepsilon'\circ_{\rmii} \varepsilon+\varepsilon'\circ_{\rmiii} \varepsilon',\\
 &&\varepsilon \circ_{\rmi} \beta+\varepsilon \circ_{\rmi} \varepsilon-\varepsilon \circ_{\rmi} \varepsilon'+\varepsilon \circ_{\rmii} \beta+\varepsilon \circ_{\rmii} \varepsilon-\varepsilon \circ_{\rmii} \varepsilon'+\varepsilon \circ_{\rmiii} \beta+\varepsilon \circ_{\rmiii} \varepsilon-\varepsilon \circ_{\rmiii} \varepsilon'.
\end{eqnarray*}
We define an isomorphism of $\mathbb{S}_2$-modules by
\begin{eqnarray*}
\eta: \antipostlie(2) \otimes \opd(2) \otimes \mathbf{k}. \sgn_{\mathbb{S}_2} &\to& \ATSu(\opd)(2)
, \quad \left\{\begin{array}{lll}
\beta \otimes \omega \otimes 1 &\mapsto& \begin{pmatrix}
  \omega\\
  \cdot
  \end{pmatrix}\smallskip
\\
\varepsilon \otimes \omega \otimes 1 &\mapsto& \begin{pmatrix}
  \omega\\
  \succ
  \end{pmatrix}
\end{array} \right .
\end{eqnarray*}
which induces an isomorphism of
$\mathbb{S}_3$-modules:
$$ \overline{\eta}: 3(\antipostlie(2) \otimes \opd(2) \otimes \mathbf{k}. \sgn_{\mathbb{S}_2})^{\otimes 2} \to 3\ATSu(\opd)^{\otimes 2}. $$
It remains to show that $\overline{\eta}$ identifies the
quadratic relation spaces of these two operads. To this end, we prove the following equations for every relation
$\gamma$ of $\opd$:
\begin{align}
\label{black-product-post-1}
& \overline{\eta}(\Psi((\beta \circ_{\rmi} \beta +\beta  \circ_{\rmii} \beta +\beta  \circ_{\rmiii} \beta)\otimes \gamma))=\ATSu_{\{x,y,z\}}(\gamma);\\
\label{black-product-post-2}
& \left\{
\begin{aligned}
& \overline{\eta}(\Psi((-\varepsilon' \circ_{\rmi} \beta+\beta \circ_{\rmiii} \varepsilon+\beta' \circ_{\rmii} \varepsilon')\otimes \gamma))=\ATSu_{\{x,y\}}(\gamma);\\
& \overline{\eta}(\Psi((-\varepsilon' \circ_{\rmii} \beta+\beta \circ_{\rmi} \varepsilon+\beta' \circ_{\rmiii} \varepsilon')\otimes \gamma))=\ATSu_{\{y,z\}}(\gamma);\\
& \overline{\eta}(\Psi((-\varepsilon' \circ_{\rmiii} \beta+\beta \circ_{\rmii} \varepsilon+\beta' \circ_{\rmi} \varepsilon')\otimes \gamma))=\ATSu_{\{x,z\}}(\gamma);
\end{aligned}
\right.\\
\label{black-product-post-3}
& \left\{
\begin{aligned}
& \overline{\eta}(\Psi((-\varepsilon \circ_{\rmii} \beta-\varepsilon \circ_{\rmii} \varepsilon+\varepsilon \circ_{\rmii} \varepsilon'-\varepsilon'\circ_{\rmiii} \varepsilon+\varepsilon'\circ_{\rmi} \varepsilon')\otimes \gamma))=\ATSu_{\{x\}}(\gamma);\\
& \overline{\eta}(\Psi((-\varepsilon \circ_{\rmiii} \beta-\varepsilon \circ_{\rmiii} \varepsilon+\varepsilon \circ_{\rmiii} \varepsilon'-\varepsilon'\circ_{\rmi} \varepsilon+\varepsilon'\circ_{\rmii} \varepsilon')\otimes \gamma))=\ATSu_{\{y\}}(\gamma);\\
& \overline{\eta}(\Psi((-\varepsilon \circ_{\rmi} \beta-\varepsilon \circ_{\rmi} \varepsilon+\varepsilon \circ_{\rmi} \varepsilon'-\varepsilon'\circ_{\rmii}
\varepsilon+\varepsilon'\circ_{\rmiii} \varepsilon')\otimes \gamma))=\ATSu_{\{z\}}(\gamma);
\end{aligned}
\right.\\
\label{black-product-post-4}
& \overline{\eta}(\Psi((\varepsilon \circ_{\rmi} \beta+\varepsilon \circ_{\rmi} \varepsilon-\varepsilon \circ_{\rmi} \varepsilon'+\varepsilon \circ_{\rmii} \beta+\varepsilon \circ_{\rmii} \varepsilon-\varepsilon \circ_{\rmii} \varepsilon'+\varepsilon \circ_{\rmiii} \beta+\varepsilon \circ_{\rmiii} \varepsilon-\varepsilon \circ_{\rmiii} \varepsilon')\otimes \gamma))\\
& \quad =\frac{1}{2} \Big(\ST(\gamma)-\sum_{\emptyset \neq J \subseteq \{x,y,z\} }  \ATSu_J (\gamma) \Big). \notag
\end{align}

By linearity, it suffices to verify
Eqs. \eqref{black-product-post-1}--\eqref{black-product-post-4}
for the monomials
$\omega\circ_{\rmi}\nu$,
$\omega\circ_{\rmii}\nu$ and
$\omega\circ_{\rmiii}\nu$,
where $\omega,\nu\in V$.
We use Eq. \eqref{eq:i-ii-iii} throughout the proof.

First note that for all $\omega$ and $\nu$ in $V$, we have
{\small
\begin{eqnarray*}
\ATSu_{\{x,y,z\}}(\omega \circ_{\rmi} \nu)=\begin{pmatrix}
  \omega\\
  \cdot
  \end{pmatrix} \circ_{\rmi} \begin{pmatrix}
  \nu\\
  \cdot
  \end{pmatrix},~
\ATSu_{\{x,y,z\}}(\omega \circ_{\rmii} \nu)=\begin{pmatrix}
  \omega\\
  \cdot
  \end{pmatrix} \circ_{\rmii} \begin{pmatrix}
  \nu\\
  \cdot
  \end{pmatrix},~
\ATSu_{\{x,y,z\}}(\omega \circ_{\rmiii} \nu)=\begin{pmatrix}
  \omega\\
  \cdot
  \end{pmatrix} \circ_{\rmiii} \begin{pmatrix}
  \nu\\
  \cdot
  \end{pmatrix}.
\end{eqnarray*}
}
Then we obtain that
\begin{eqnarray*}
\overline{\eta}(\Psi((\beta \circ_{\rmi} \beta +\beta  \circ_{\rmii} \beta +\beta  \circ_{\rmiii} \beta)\otimes (\omega \circ_{\rmi} \nu)))=\ATSu_{\{x,y,z\}}(\omega \circ_{\rmi} \nu).
\end{eqnarray*}
Second from the following equalities
{\small
\begin{eqnarray*}
\ATSu_{\{x,y\}}(\omega \circ_{\rmi} \nu)=\begin{pmatrix}
  \omega\\
  \prec
  \end{pmatrix} \circ_{\rmi} \begin{pmatrix}
  \nu\\
  \cdot
  \end{pmatrix},~
\ATSu_{\{x,y\}}(\omega \circ_{\rmii} \nu)=\begin{pmatrix}
  \omega\\
  \cdot
  \end{pmatrix} \circ_{\rmii} \begin{pmatrix}
  \nu\\
  \prec
  \end{pmatrix},~
\ATSu_{\{x,y\}}(\omega \circ_{\rmiii} \nu)=\begin{pmatrix}
  \omega\\
  \cdot
  \end{pmatrix} \circ_{\rmiii} \begin{pmatrix}
  \nu\\
  \succ
  \end{pmatrix},
\end{eqnarray*}
}
we also have
\begin{eqnarray*}
\overline{\eta}(\Psi((-\varepsilon' \circ_{\rmi} \beta+\beta \circ_{\rmiii} \varepsilon+\beta' \circ_{\rmii} \varepsilon')\otimes (\omega \circ_{\rmi} \nu)))=\ATSu_{\{x,y\}}(\omega \circ_{\rmi} \nu).
\end{eqnarray*}
Third, since we have
{\small
\begin{eqnarray*}
\ATSu_{\{x\}}(\omega \circ_{\rmi} \nu)=\begin{pmatrix}
  \omega\\
  \prec
  \end{pmatrix} \circ_{\rmi} \begin{pmatrix}
  \nu\\
  \prec
  \end{pmatrix},~
\ATSu_{\{x\}}(\omega \circ_{\rmii} \nu)=\begin{pmatrix}
  \omega\\
  \succ
  \end{pmatrix} \circ_{\rmii} \begin{pmatrix}
  \nu\\
  \star
  \end{pmatrix},~
\ATSu_{\{x\}}(\omega \circ_{\rmiii} \nu)=\begin{pmatrix}
  \omega\\
  \prec
  \end{pmatrix} \circ_{\rmiii} \begin{pmatrix}
  \nu\\
  \succ
  \end{pmatrix},
\end{eqnarray*}
}
which implies that
\begin{eqnarray*}
\overline{\eta}(\Psi((-\varepsilon \circ_{\rmii} \beta-\varepsilon \circ_{\rmii} \varepsilon+\varepsilon \circ_{\rmii} \varepsilon'-\varepsilon'\circ_{\rmiii} \varepsilon+\varepsilon'\circ_{\rmi} \varepsilon')\otimes (\omega \circ_{\rmi} \nu)))=
\ATSu_{\{x\}}(\omega \circ_{\rmi} \nu).
\end{eqnarray*}
Finally, we also have
\[
\begin{aligned}
\sum_{\emptyset\neq J\subseteq\{x,y,z\}}
\ATSu_J(\omega\circ_{\rmi}\nu)
&=\left(
\begin{pmatrix}\omega\\ \prec\end{pmatrix}
+
\begin{pmatrix}\omega\\ \cdot\end{pmatrix}
\right)
\circ_{\rmi}
\left(
\begin{pmatrix}\nu\\ \prec\end{pmatrix}
+
\begin{pmatrix}\nu\\ \succ\end{pmatrix}
+
\begin{pmatrix}\nu\\ \cdot\end{pmatrix}
\right)
+
\begin{pmatrix}\omega\\ \succ\end{pmatrix}
\circ_{\rmi}
\begin{pmatrix}\nu\\ \star\end{pmatrix}
\\
&=
\ST(\omega\circ_{\rmi}\nu)
+
2\begin{pmatrix}\omega\\ \succ\end{pmatrix}
\circ_{\rmi}
\begin{pmatrix}\nu\\ \star\end{pmatrix},
\end{aligned}
\]
which implies that
\[
\begin{aligned}
&\overline{\eta}\Bigl(
\Psi\Bigl(
\Bigl(
\varepsilon\circ_{\rmi}\beta
+\varepsilon\circ_{\rmi}\varepsilon
-\varepsilon\circ_{\rmi}\varepsilon'
+\varepsilon\circ_{\rmii}\beta
+\varepsilon\circ_{\rmii}\varepsilon
-\varepsilon\circ_{\rmii}\varepsilon'
+\varepsilon\circ_{\rmiii}\beta
+\varepsilon\circ_{\rmiii}\varepsilon
-\varepsilon\circ_{\rmiii}\varepsilon'
\Bigr)\otimes(\omega\circ_{\rmi}\nu)
\Bigr)
\Bigr)
\\
&=
\begin{pmatrix}\omega\\ \succ\end{pmatrix}
\circ_{\rmi}
\left(
\begin{pmatrix}\nu\\ \cdot\end{pmatrix}
+
\begin{pmatrix}\nu\\ \succ\end{pmatrix}
+
\begin{pmatrix}\nu\\ \prec\end{pmatrix}
\right)
 =
\frac12\left(
\ST(\omega\circ_{\rmi}\nu)
-
\sum_{\emptyset\neq J\subseteq\{x,y,z\}}
\ATSu_J(\omega\circ_{\rmi}\nu)
\right).
\end{aligned}
\]
The corresponding identities for $\omega\circ_{\rmii}\nu$ and $\omega\circ_{\rmiii}\nu$ follow by analogous calculations in each of the four cases. Since these monomials span $\huaT(V)(3)$, Eqs.~\eqref{black-product-post-1}--\eqref{black-product-post-4}
hold for every $\gamma\in\huaT(V)(3)$.

By Eqs.~\eqref{black-product-post-1}--\eqref{black-product-post-4},
the eight images span the same space as
\[
\bigl\{
\ATSu_J(\gamma)
\mid
\emptyset\neq J\subseteq\{x,y,z\}
\bigr\}
\cup
\{\ST(\gamma)\}.
\]
Letting $\gamma$ range over the relation space of $\opd$,
we deduce that $\overline{\eta}$ identifies the two quadratic
relation spaces. Hence $\eta$ induces the required isomorphism
of operads.
\end{proof}

We have the following counterpart of Theorem \ref{thm:Manin-product-post} for anti-successors.

\begin{thm}\label{thm:Manin-product}
Let $\opd$ be a binary quadratic operad. There is an isomorphism of operads
$$ \ABSu(\opd) \cong \antiprelie \bullet \opd. $$
\end{thm}

\begin{proof}
The proof is analogous to that of Theorem~\ref{thm:Manin-product-post};
we omit the details.
\end{proof}

In the nonsymmetric setting, we obtain the following analogue of Theorems \ref{thm:Manin-product-post} and \ref{thm:Manin-product}, which can be proved by similar arguments.

\begin{pro}
Let $\opd$ be a binary quadratic nonsymmetric operad. There are isomorphisms of nonsymmetric operads
$$ \ABSu(\opd) \cong \antiden~\blacksquare~\opd, \quad \ATSu(\opd) \cong \antitriden~\blacksquare~\opd. $$
Here $\blacksquare$ denotes the black square product in \cite{FG,Val2}.
\end{pro}

Theorems \ref{thm:Manin-product-post} and \ref{thm:Manin-product} give an efficient way to compute the Manin black product of a binary operad with the operad $\antiprelie$ or $\antipostlie$, as illustrated by the following corollaries.

\begin{cor}
We have $\antiprelie \bullet \com\cong \antizinbiel$, $\antiprelie \bullet \ass\cong \antiden$ and  $\antiprelie \bullet \poisson\cong \antiprepoisson$.
\end{cor}

\begin{cor}
We have $\antipostlie \bullet \com\cong \antipostcom$, $\antipostlie \bullet \lie\cong \antipostlie$ and $\antitriden~\blacksquare~\ass\cong\antitriden$.
\end{cor}

\section{The anti-replicators of a binary operad}\label{sec:anti-rep}

We first introduce the notions of the
canonical quadratic anti-duplicator and anti-triplicator of a
binary quadratic operad, collectively referred to as the {\bf
canonical quadratic anti-replicators}. We then consider their
combinatorial counterparts for binary quadratic operads, defined
via the anti-duplicator and anti-triplicator constructions on
labeled planar binary trees. We compare these two constructions
and determine when they agree. Finally, we apply the combinatorial construction to obtain anti-dupli-com algebras and anti-dupli-Lie algebras, and use the operad $\mathrm{2Ass}$ to show that the two constructions differ in general.

\subsection{The canonical quadratic anti-replicators of a binary operad}

Let $\opd=\mathcal{T}(\gensp)/ (\relsp)$ be a binary (symmetric) quadratic operad with a finite-dimensional $\mathbb{S}_2$-module $V$ and let $R\subseteq \mathcal{T}(\gensp)(3)$ be an $\mathbb{S}_3$-submodule. We define the Koszul dual of $V$ by $V^{\vee}=V^{\ast} \otimes \sgn_2$. For a basis $\{e_1,\cdots,e_n\}$ of $V$ and its dual basis $\{e_1^{\vee},\cdots,e_n^{\vee}\}$ of $V^{\vee}$, there is a natural quadratic pairing with respect to the duality given by
\begin{eqnarray*}
&&\langle\cdot,\cdot\rangle: \mathcal{T}(\gensp^{\vee})(3) \otimes  \mathcal{T}(\gensp)(3) \to \K,\\
&& \langle e_i^{\vee}\circ_{u} e_j^{\vee}, e_k \circ_{v} e_l \rangle=\delta_{(i,k)}\delta_{(j,l)}\delta_{(u,v)} \in \K,\quad \forall u,v\in \{\rmi,\rmii,\rmiii\}.
\end{eqnarray*}

Let $\gensp$ be a vector space with a basis $\genbas$.
We introduce two vector spaces
\begin{eqnarray}\label{gendup}
\gendup
&=&
\gensp\otimes
(\mathbf{k}\dashv\oplus\mathbf{k}\vdash),\\
\label{gentrip}
\gentrip
&=&
\gensp\otimes
(\mathbf{k}\dashv\oplus\mathbf{k}\vdash
\oplus\mathbf{k}\perp).
\end{eqnarray}
For $\omega\in\gensp$, we denote
$\omega~\otimes\dashv$, $\omega~\otimes\vdash$
and $\omega~\otimes\perp$ by
$
\begin{pmatrix}\omega\\ \dashv\end{pmatrix},
~
\begin{pmatrix}\omega\\ \vdash\end{pmatrix},
~
\begin{pmatrix}\omega\\ \perp\end{pmatrix},
$
respectively. Consequently,
$
\left\{
\begin{pmatrix}\omega\\ \dashv\end{pmatrix},
\begin{pmatrix}\omega\\ \vdash\end{pmatrix}
\,\middle|\,\omega\in\genbas
\right\}
$
and
$
\left\{
\begin{pmatrix}\omega\\ \dashv\end{pmatrix},
\begin{pmatrix}\omega\\ \vdash\end{pmatrix},
\begin{pmatrix}\omega\\ \perp\end{pmatrix}
\,\middle|\,\omega\in\genbas
\right\}
$
form bases of $\gendup$ and $\gentrip$,
respectively. Both spaces are regarded
as arity-graded vector spaces concentrated
in arity $2$.

\begin{defi}\label{def-can-rep}
Retain the above notations.
\begin{enumerate}
\item[{\rm(i)}] For $x\in \{1,2,3\}$, let
  $$ q_x^V: \mathcal{T}(\gendup)(3) \to \mathcal{T}(\gensp)(3),\quad t^V: \mathcal{T}(\gendup)(3) \to \mathcal{T}(\gensp)(3) $$
  be the adjoints of $\ABSu_x$ and $\st$ acting on $\mathcal{T}(\gensp^{\vee})(3)$, respectively. That means, for all $h\in \mathcal{T}(\gendup)(3)$ and $r\in \mathcal{T}(\gensp^{\vee})(3)$, we have
  \begin{eqnarray}\label{eq:adjunction-binary-revised}
  \langle q_x^V(h),r\rangle=\langle h, \ABSu_x(r) \rangle,\quad \langle t^V(h),r\rangle=\langle h, \st(r) \rangle.
  \end{eqnarray}
  The {\bf canonical quadratic anti-duplicator} of $\opd$ is defined to be the binary operad
  $$ \widehat{\ADup}(\opd):=\mathcal{T}(\gendup)/ (L_2^V(R)), $$
  where the quadratic relation subspace is given by
  $$ L_2^V (R):= \bigcap_{x=1}^3 (q_x^V)^{-1} (R) \cap (t^V)^{-1} (R) \subseteq  \mathcal{T}(\gendup)(3). $$
  \item[{\rm(ii)}] For each nonempty subset $J \subseteq \{1,2,3\}$, let
  $$ q_J^V: \mathcal{T}(\gentrip)(3) \to \mathcal{T}(\gensp)(3),\quad T^V: \mathcal{T}(\gentrip)(3) \to \mathcal{T}(\gensp)(3) $$
  be the adjoints of $\ATSu_J$ and $\ST$ acting on $\mathcal{T}(\gensp^{\vee})(3)$, respectively. That means, for all $h\in \mathcal{T}(\gentrip)(3)$ and $r\in \mathcal{T}(\gensp^{\vee})(3)$, we have
  \begin{eqnarray}\label{eq:adjunction-ternary-revised}
  \langle q_J^V(h),r\rangle=\langle h, \ATSu_J (r) \rangle,\quad \langle T^V(h),r\rangle=\langle h, \ST(r) \rangle.
  \end{eqnarray}
  The {\bf canonical quadratic anti-triplicator} of $\opd$ is defined to be the binary operad
  $$ \widehat{\ATrip}(\opd):=\mathcal{T}(\gentrip)/ (L_3^V(R)), $$
  where the quadratic relation subspace is given by
  $$ L_3^V (R):= \bigcap_{\emptyset \neq J \subseteq \{1,2,3\}}  (q_J^V)^{-1} (R) \cap (T^V)^{-1} (R) \subseteq  \mathcal{T}(\gentrip)(3). $$
\end{enumerate}
Here $(L_m^V(R))$ denotes the operadic ideal generated by the relation subspace $L_m^V(R)$, for $m=2,3$.
\end{defi}

\begin{rmk}
The maps $q_x^V$, $t^V$, $q_J^V$ and $T^V$ in Definition \ref{def-can-rep} are defined by adjunction and do not require a basis of $V$. The natural identifications $(\widetilde{V})^\vee\cong \overline{V^\vee}$ and  $(\widehat{V})^\vee\cong \overrightarrow{V^\vee}$ are used here. Since $R\subseteq \mathcal{T}(\gensp)(3)$ is an $\mathbb{S}_3$-submodule, it is easy to check that $L_2^V (R)$ and $L_3^V (R)$ are also $\mathbb{S}_3$-submodules.
\end{rmk}

\subsection{The combinatorial anti-replicators of a binary operad}

\subsubsection{Anti-duplicators of labeled planar binary trees}

\begin{defi}\label{def:anti-dup}{\rm
Let $\gensp$ be a vector space with a basis $\genbas$. For a labeled $n$-tree $\tau$ in $\calt(\genbas)$, define $\overline{\tau}$ in $\mathcal{T}_{\hspace*{-0.1cm}ns}(\gendup)$ as follows:
\begin{itemize}
 \item[$\bullet$] $  \overline{~|~}=|$
 \item[$\bullet$] when $n\geq 2$, $\overline{\tau}$ is obtained by replacing each decoration $\gop \in \vin(\tau)$ by
     $$\begin{pmatrix}
  \gop \\
  \dagger
  \end{pmatrix}:=\Biggl\{-\begin{pmatrix}
  \gop \\
  \dashv
  \end{pmatrix},-\begin{pmatrix}
  \gop \\
  \vdash
  \end{pmatrix}\Biggr\}.$$
\end{itemize}
Thus, $\overline{\tau}$ is a set of labeled trees. We then extend this definition to $\mathcal{T}_{\hspace*{-0.1cm}ns}(\gendup)$ by linearity.
}\end{defi}

\begin{defi}
Let $\gensp$ be a vector space with a basis $\genbas$. Let $\tau$ be a labeled $n$-tree in
$\calt(\genbas)$. The {\bf anti-duplicator} $\ADup_{x}(\tau)$ of $\tau$ with respect to a leaf $x\in \lin(\tau)$ is a subset of $\mathcal{T}_{\hspace*{-0.1cm}ns}(\gendup)$ defined by induction on
$|\lin(\tau)|$ as follows:
\begin{itemize}
 \item[$\bullet$] $\ADup_{x}(~|~)=\{~|~\}$~;
 \item[$\bullet$] Assume that $\ADup_{x}(\tau)$ have been defined for $\tau$ with $|\lin(\tau)|\leq k$ for any $k\geq 1$.  Then for a labeled $(k+1)$-tree $\tau\in \calt(\genbas)$ with its decomposition $\tau_{l} \vee_{\gop} \tau_{r}$, we define
     $$
     \ADup_{x}(\tau)= \ADup_{x}(\tau_{l} \vee_{\gop} \tau_{r}) = \left\{\begin{array}{ll}  \ADup_{x}(\tau_{l}) \vee_{\begin{pmatrix}
  \gop \\
  \dashv
  \end{pmatrix}} \overline{\tau_{r}},\quad & \text{if}~x \in \lin(\tau_{l});  \\
     \overline{\tau_{l}} \vee_{\begin{pmatrix}
  \gop \\
  \vdash
  \end{pmatrix}} \ADup_{x}(\tau_{r}), \quad & \text{if}~x \in \lin(\tau_{r}). \end{array}\right.
     $$
\end{itemize}
\end{defi}

For labeled $n$-trees $\tau_i$ with the same set of leaf decorations and $c_i\in \K$ for all $1\leq i \leq r$, define
\begin{eqnarray}
\ADup_{x} (\sum\limits_{i=1}^r c_i \tau_i):=\sum\limits_{i=1}^r c_i \ADup_{x} (\tau_i).
\end{eqnarray}
Throughout the remainder of the paper, we use the following notation:
\begin{eqnarray}\label{eq:convention}
\sum\limits_{i=1}^r c_i W_i:=\Biggl\{ \sum\limits_{i=1}^r c_i w_i~\Biggl|~w_i\in W_i,\quad 1\leq i \leq r  \Biggr\},
\end{eqnarray}
for nonempty subsets $W_i$ of a $\K$-module.

The following result characterizes the
anti-duplicators of labeled trees explicitly.

\begin{pro}\label{pro:ADup}
Let $\gensp$ be a vector space with a basis $\genbas$, $\tau$ be a
labeled $n$-tree in $\calt(\genbas)$ and $x$ be in $\lin(\tau)$.
The anti-duplicator $\ADup_{x}(\tau)$ of $\tau$ is obtained by
relabelling each vertex of $\tau$ according to the following
rules:
\begin{enumerate}
 \item[\rm{(a)}]
We replace the label $\gop$ of each vertex on the path from the root to the leaf $x$ of $\tau$ by
\begin{itemize}
 \item[\rm{(i)}] $\begin{pmatrix}
  \gop \\
  \dashv
  \end{pmatrix}$,\quad if the path turns left at this vertex;
\item[\rm{(ii)}] $\begin{pmatrix}
  \gop \\
  \vdash
  \end{pmatrix}$,\quad if the path turns right at this vertex;
\end{itemize}
 \item[\rm{(b)}]
We replace the label $\gop$ of each vertex not on the path from the root to the leaf $x$ of $\tau$ by
$$\begin{pmatrix}
  \gop \\
  \dagger
  \end{pmatrix}:=\Biggl\{-\begin{pmatrix}
  \gop \\
  \dashv
  \end{pmatrix},-\begin{pmatrix}
  \gop \\
  \vdash
  \end{pmatrix}\Biggr\}.$$
\end{enumerate}
\end{pro}

\begin{proof}
We can prove it by induction on $|\lin(\tau)|\geq 1$. We omit the details for brevity.
\end{proof}

\begin{ex}
{\tiny$$\ADup_{x_{2}} \left( \begin{xy}
(0,8)*+{x_{1}}="v1";(10,8)*+{x_{2}}="v2";(20,8)*+{x_{3}}="v3"; (30,8)*+{x_{4}}="v4";
(5,0)*+{\omega_{1}}="v5";(25,0)*+{\omega_{2}}="v6"; (15,-6)*+{\omega_{3}}="v7";(15,-12)*+{}="v8";
{\ar@{-}"v5";"v1"};{\ar@{-}"v5";"v2"};{\ar@{-}"v6";"v3"};{\ar@{-}"v6";"v4"};{\ar@{-}"v7";"v5"};{\ar@{-}"v7";"v6"};{\ar@{-}"v8";"v7"}
\end{xy} \right)= \begin{xy}
(0,10)*+{x_{1}}="v1";(10,10)*+{x_{2}}="v2";(20,10)*+{x_{3}}="v3"; (30,10)*+{x_{4}}="v4";
(5,0)*+{\svec{\omega_{1}}{\vdash}}="v5";(25,0)* +{\svec{\omega_{2}}{\dagger}}="v6"; (15,-6)*+{\svec{\omega_{3}}{\dashv}}="v7";(15,-15)*+{}="v8";
{\ar@{-}"v5";"v1"};{\ar@{->}"v5";"v2"};{\ar@{-}"v6";"v3"};{\ar@{-}"v6";"v4"};{\ar@{->}"v7";"v5"};{\ar@{-}"v7";"v6"};{\ar@{->}"v8";"v7"}
\end{xy}
= \left\{-\begin{xy}
(0,10)*+{x_{1}}="v1";(10,10)*+{x_{2}}="v2";(20,10)*+{x_{3}}="v3"; (30,10)*+{x_{4}}="v4";
(5,0)*+{\svec{\omega_{1}}{\vdash}}="v5";(25,0)* +{\svec{\omega_{2}}{\dashv}}="v6"; (15,-6)*+{\svec{\omega_{3}}{\dashv}}="v7";(15,-15)*+{}="v8";
{\ar@{-}"v5";"v1"};{\ar@{-}"v5";"v2"};{\ar@{-}"v6";"v3"};{\ar@{-}"v6";"v4"};{\ar@{-}"v7";"v5"};{\ar@{-}"v7";"v6"};{\ar@{-}"v8";"v7"}\end{xy},\quad
-\begin{xy}
(0,10)*+{x_{1}}="v1";(10,10)*+{x_{2}}="v2";(20,10)*+{x_{3}}="v3"; (30,10)*+{x_{4}}="v4";
(5,0)*+{\svec{\omega_{1}}{\vdash}}="v5";(25,0)* +{\svec{\omega_{2}}{\vdash}}="v6"; (15,-6)*+{\svec{\omega_{3}}{\dashv}}="v7";(15,-15)*+{}="v8";
{\ar@{-}"v5";"v1"};{\ar@{-}"v5";"v2"};{\ar@{-}"v6";"v3"};{\ar@{-}"v6";"v4"};{\ar@{-}"v7";"v5"};{\ar@{-}"v7";"v6"};{\ar@{-}"v8";"v7"}
\end{xy}\right\}
$$}
\end{ex}

\subsubsection{Anti-triplicators of labeled planar binary trees}

\begin{defi}\label{def:anti-trip}{\rm
Let $\gensp$ be a vector space with a basis $\genbas$. For a labeled $n$-tree $\tau$ in $\calt(\genbas)$, define $\overrightarrow{\tau}$ in $\mathcal{T}_{\hspace*{-0.1cm}ns}(\gentrip)$ as follows:
\begin{itemize}
 \item[$\bullet$] $  \overrightarrow{~|~}=|$
 \item[$\bullet$] when $n\geq 2$, $\overrightarrow{\tau}$ is obtained by replacing each decoration $\gop \in \vin(\tau)$ by
     $$\begin{pmatrix}
  \gop \\
  \ddagger
  \end{pmatrix}:=\Biggl\{-\begin{pmatrix}
  \gop \\
  \dashv
  \end{pmatrix},-\begin{pmatrix}
  \gop \\
  \vdash
  \end{pmatrix},-\begin{pmatrix}
  \gop \\
  \perp
  \end{pmatrix}\Biggr\}.$$
\end{itemize}
Thus, $\overrightarrow{\tau}$ is a set of labeled trees. We then extend this definition to $\mathcal{T}_{\hspace*{-0.1cm}ns}(\gentrip)$ by linearity.
}\end{defi}

\begin{defi}
Let $\gensp$ be a vector space with a basis $\genbas$. Let $\tau$ be a labeled $n$-tree in
$\calt(\genbas)$. The {\bf anti-triplicator} $\ATrip_{J}(\tau)$ of $\tau$ with respect to a subset $J \subseteq \lin(\tau)$ is a subset of $\mathcal{T}_{\hspace*{-0.1cm}ns}(\gentrip)$ defined by induction on $|\lin(\tau)|$ as follows:
\begin{itemize}
 \item[$\bullet$] $\ATrip_{J}(~|~)=\{~|~\}$~;
 \item[$\bullet$] $\ATrip_{\emptyset}(\tau)=\overrightarrow{\tau}$;
 \item[$\bullet$] Assume that $\ATrip_{J}(\tau)$ have been defined for $\tau$ with $|\lin(\tau)|\leq k$ for any $k\geq 1$.  Then for a labeled $(k+1)$-tree $\tau\in \calt(\genbas)$ with its decomposition $\tau_{l} \vee_{\gop} \tau_{r}$, we define
     $$
     \ATrip_{J}(\tau)= \ATrip_{J}(\tau_{l} \vee_{\gop} \tau_{r}) = \left\{\begin{array}{ll}  \ATrip_{J}(\tau_{l}) \vee_{\begin{pmatrix}
  \gop \\
  \dashv
  \end{pmatrix}} \overrightarrow{\tau_{r}},\quad & \text{if}~J \subseteq \lin(\tau_l);  \\
     \overrightarrow{\tau_{l}} \vee_{\begin{pmatrix}
  \gop \\
  \vdash
  \end{pmatrix}} \ATrip_{J}(\tau_{r}), \quad & \text{if}~J \subseteq \lin(\tau_r);\\
  \overrightarrow{\tau_{l}} \vee_{\begin{pmatrix}
  \gop \\
  \ddagger
  \end{pmatrix}} \overrightarrow{\tau_{r}},\quad & \text{if}~J=\emptyset;  \\
  \ATrip_{J\cap \lin(\tau_l)}(\tau_{l}) \vee_{\begin{pmatrix}
  \gop \\
  \perp
  \end{pmatrix}} \ATrip_{J\cap \lin(\tau_r)}(\tau_{r}),\quad & \text{otherwise.}
  \end{array}\right.
     $$
\end{itemize}
\end{defi}

The following result gives an explicit
characterization of the anti-triplicators.

\begin{pro}
Let $\gensp$ be a vector space with a basis $\genbas$, $\tau$ be a
labeled $n$-tree in $\calt(\genbas)$ and $J$ be a nonempty subset
of $\lin(\tau)$. The anti-triplicator $\ATrip_{J}(\tau)$ of $\tau$
is obtained by relabelling each vertex of $\vin(\tau)$ by the
following rules:
\begin{enumerate}
 \item[\rm{(a)}]
We replace the label $\gop$ of each vertex on the paths from the root of $\tau$ to some leaves (possibly multiple) in $J$ by
\begin{itemize}
 \item[\rm{(i)}] $\begin{pmatrix}
  \gop \\
  \dashv
  \end{pmatrix}$,\quad if all such paths turn left at this vertex;
\item[\rm{(ii)}] $\begin{pmatrix}
  \gop \\
  \vdash
  \end{pmatrix}$,\quad if all such paths turn right at this vertex;
\item[\rm{(iii)}] $\begin{pmatrix}
  \gop \\
  \perp
  \end{pmatrix}$,\quad if some of these paths turn left and some of these paths turn right.
\end{itemize}
 \item[\rm{(b)}]
We replace the label $\gop$ of each vertex not on the path from the root to any leaf $x\in J$ by
$$\begin{pmatrix}
  \gop \\
  \ddagger
  \end{pmatrix}:=\Biggl\{-\begin{pmatrix}
  \gop \\
  \dashv
  \end{pmatrix},-\begin{pmatrix}
  \gop \\
  \vdash
  \end{pmatrix},-\begin{pmatrix}
  \gop \\
  \perp
  \end{pmatrix}\Biggr\}.$$
\end{enumerate}
\end{pro}

\begin{proof}
It follows from an induction on $|\lin(\tau)|\geq 1$.
\end{proof}

\begin{ex}
{\footnotesize$$\ATrip_{\{x_{1},x_{2}\}} \left( \begin{xy}
(0,8)*+{x_{1}}="v1";(10,8)*+{x_{2}}="v2";(20,8)*+{x_{3}}="v3"; (30,8)*+{x_{4}}="v4";
(5,0)*+{\omega_{1}}="v5";(25,0)*+{\omega_{2}}="v6"; (15,-6)*+{\omega_{3}}="v7";(15,-12)*+{}="v8";
{\ar@{-}"v5";"v1"};{\ar@{-}"v5";"v2"};{\ar@{-}"v6";"v3"};{\ar@{-}"v6";"v4"};{\ar@{-}"v7";"v5"};{\ar@{-}"v7";"v6"};{\ar@{-}"v8";"v7"}
\end{xy} \right)= \begin{xy}
(0,10)*+{x_{1}}="v1";(10,10)*+{x_{2}}="v2";(20,10)*+{x_{3}}="v3"; (30,10)*+{x_{4}}="v4";
(5,0)*+{\svec{\omega_{1}}{\perp}}="v5";(25,0)* +{\svec{\omega_{2}}{\ddagger}}="v6"; (15,-6)*+{\svec{\omega_{3}}{\dashv}}="v7";(15,-15)*+{}="v8";
{\ar@{->}"v5";"v1"};{\ar@{->}"v5";"v2"};{\ar@{-}"v6";"v3"};{\ar@{-}"v6";"v4"};{\ar@{->}"v7";"v5"};{\ar@{-}"v7";"v6"};{\ar@{->}"v8";"v7"}
\end{xy}
$$
$$
= \left\{-\begin{xy}
(0,10)*+{x_{1}}="v1";(10,10)*+{x_{2}}="v2";(20,10)*+{x_{3}}="v3"; (30,10)*+{x_{4}}="v4";
(5,0)*+{\svec{\omega_{1}}{\perp}}="v5";(25,0)* +{\svec{\omega_{2}}{\dashv}}="v6"; (15,-6)*+{\svec{\omega_{3}}{\dashv}}="v7";(15,-15)*+{}="v8";
{\ar@{-}"v5";"v1"};{\ar@{-}"v5";"v2"};{\ar@{-}"v6";"v3"};{\ar@{-}"v6";"v4"};{\ar@{-}"v7";"v5"};{\ar@{-}"v7";"v6"};{\ar@{-}"v8";"v7"}\end{xy},\quad
-\begin{xy}
(0,10)*+{x_{1}}="v1";(10,10)*+{x_{2}}="v2";(20,10)*+{x_{3}}="v3"; (30,10)*+{x_{4}}="v4";
(5,0)*+{\svec{\omega_{1}}{\perp}}="v5";(25,0)* +{\svec{\omega_{2}}{\vdash}}="v6"; (15,-6)*+{\svec{\omega_{3}}{\dashv}}="v7";(15,-15)*+{}="v8";
{\ar@{-}"v5";"v1"};{\ar@{-}"v5";"v2"};{\ar@{-}"v6";"v3"};{\ar@{-}"v6";"v4"};{\ar@{-}"v7";"v5"};{\ar@{-}"v7";"v6"};{\ar@{-}"v8";"v7"}
\end{xy},\quad
-\begin{xy}
(0,10)*+{x_{1}}="v1";(10,10)*+{x_{2}}="v2";(20,10)*+{x_{3}}="v3"; (30,10)*+{x_{4}}="v4";
(5,0)*+{\svec{\omega_{1}}{\perp}}="v5";(25,0)* +{\svec{\omega_{2}}{\perp}}="v6"; (15,-6)*+{\svec{\omega_{3}}{\dashv}}="v7";(15,-15)*+{}="v8";
{\ar@{-}"v5";"v1"};{\ar@{-}"v5";"v2"};{\ar@{-}"v6";"v3"};{\ar@{-}"v6";"v4"};{\ar@{-}"v7";"v5"};{\ar@{-}"v7";"v6"};{\ar@{-}"v8";"v7"}
\end{xy}\right\}
$$
}
\end{ex}

\subsubsection{The combinatorial anti-replicators of a binary operad}

\begin{defi}\label{def:ADup-sym-operad}
Retain the conditions in Definition~\ref{def:ATBSu-sym-operad}.
\begin{enumerate}
\item[{\rm(i)}]
The {\bf combinatorial anti-duplicator} of $\opd$ is defined to be the binary quadratic operad
$$
\ADup(\opd):=\mathcal{T}(\gendup)/
\bigl( \ADup(\relsp)\cap \comp(\relsp)\bigr),
$$
where $\gendup=\gensp\otimes(\mathbf{k}\dashv\oplus\mathbf{k}\vdash)$ is concentrated in arity $2$, and the $\mathbb{S}_2$-action on $\gendup$ is given by
$$
\begin{pmatrix}\gop\\ \dashv\end{pmatrix}^{(12)}
:=\begin{pmatrix}\gop^{(12)}\\ \vdash\end{pmatrix},
\quad
\begin{pmatrix}\gop\\ \vdash\end{pmatrix}^{(12)}
:=\begin{pmatrix}\gop^{(12)}\\ \dashv\end{pmatrix},
\quad \forall\,\gop\in\gensp.
$$
The space of relations is generated as an $\mathbb{S}$-module by
$$
\ADup(\relsp):=\spann_{\mathbf{k} [\mathbb{S}_3]}
\left( \bigcup_{s=1}^k\bigcup_{x\in\lin(r_s)}
\ADup_x(r_s) \right),~\text{where }\ADup_x(r_s):=\sum_i c_{s,i}\ADup_x(\tau_{s,i}),
$$
and
$$
\comp(\relsp):=\spann_{\mathbf{k} [\mathbb{S}_3]} \left( \bigcup_{s=1}^k\comp(r_s) \right),~
\text{where }\comp(r_s):=\sum_i c_{s,i}\overline{\tau_{s,i}},
$$
according to the notation by Eq.~\eqref{eq:convention}. Here $\comp(\tau_{s,i})=\overline{\tau_{s,i}}$ is obtained by replacing each vertex decoration $\gop$ of $\tau_{s,i}$ by
$$
\begin{pmatrix}\gop\\ \dagger\end{pmatrix}
:=\Biggl\{
-\begin{pmatrix}\gop\\ \dashv\end{pmatrix},
-\begin{pmatrix}\gop\\ \vdash\end{pmatrix}
\Biggr\}.
$$

\item[{\rm(ii)}]
The {\bf combinatorial anti-triplicator} of $\opd$ is defined to be the binary quadratic operad
$$
\ATrip(\opd):=\mathcal{T}(\gentrip)/
\bigl(\ATrip(\relsp)\cap \COMP(\relsp)\bigr),
$$
where $\gentrip=\gensp\otimes
(\mathbf{k}\dashv\oplus\mathbf{k}\vdash\oplus\mathbf{k}\perp)$
is concentrated in arity $2$, and the $\mathbb{S}_2$-action on $\gentrip$ is given by
$$
\begin{pmatrix}\gop\\ \dashv\end{pmatrix}^{(12)}
:=\begin{pmatrix}\gop^{(12)}\\ \vdash\end{pmatrix},
\quad
\begin{pmatrix}\gop\\ \vdash\end{pmatrix}^{(12)}
:=\begin{pmatrix}\gop^{(12)}\\ \dashv\end{pmatrix},
\quad
\begin{pmatrix}\gop\\ \perp\end{pmatrix}^{(12)}
:=\begin{pmatrix}\gop^{(12)}\\ \perp\end{pmatrix},
\quad \forall\,\gop\in\gensp.
$$
The space of relations is generated as an $\mathbb{S}$-module by
$$
\ATrip(\relsp):=\spann_{\mathbf{k} [\mathbb{S}_3]}
\left( \bigcup_{s=1}^k
\bigcup_{\emptyset\neq J\subseteq\lin(r_s)}
\ATrip_J(r_s)  \right),~\text{where }\ATrip_J(r_s):=\sum_i c_{s,i}\ATrip_J(\tau_{s,i}),
$$
and
$$
\COMP(\relsp):=\spann_{\mathbf{k} [\mathbb{S}_3]}
\left(\bigcup_{s=1}^k\COMP(r_s) \right),~\text{where } \COMP(r_s):=\sum_i c_{s,i}\overrightarrow{\tau_{s,i}}.
$$
Here $\COMP(\tau_{s,i})=\overrightarrow{\tau_{s,i}}$ is obtained by replacing each vertex decoration $\gop$ of $\tau_{s,i}$ by
$$
\begin{pmatrix}\gop\\ \ddagger\end{pmatrix}
:=\Biggl\{
-\begin{pmatrix}\gop\\ \dashv\end{pmatrix},
-\begin{pmatrix}\gop\\ \vdash\end{pmatrix},
-\begin{pmatrix}\gop\\ \perp\end{pmatrix}
\Biggr\}.
$$
\end{enumerate}
\end{defi}

\begin{rmk}
The above definitions also apply to binary nonsymmetric operads by omitting the symmetric group actions and replacing the free symmetric operad $\mathcal{T}(\gensp)$ with the free nonsymmetric operad $\mathcal{T}_{\hspace*{-0.1cm}ns}(\gensp)$.
\end{rmk}

The following theorem compares the combinatorial and canonical quadratic anti-replicators. In particular, it shows that these two constructions coincide under the full-support condition.

Let $V$ be a finite-dimensional $\mathbb{S}_2$-module. A basis $\mathcal B$ of $V$ is said to be
{\bf compatible with the $\mathbb{S}_2$-action}, if for every $b\in\mathcal B$, there exist
$b'\in\mathcal B$ and
$\varepsilon_b\in\{1,-1\}$ such that $b^{(12)}=\varepsilon_b b'$. In other words, the action of $\mathbb{S}_2$
permutes the basis elements up to signs.
Such a basis induces a monomial basis of
$\mathcal{T}(V)(3)$.

\begin{thm}\label{lem:comparison-two-rep}
Let $V$ be a finite-dimensional $\mathbb{S}_2$-module
with a compatible basis $\mathcal B$, and let
$R\subseteq\huaT(V)(3)$ be an $\mathbb{S}_3$-submodule.
Denote by $\{e_a\}$ the corresponding monomial basis
of $\huaT(V)(3)$ induced by $\mathcal B$, and by
$\operatorname{supp}_{\mathcal B}(R)$ the set of monomials
that occur with nonzero coefficients in at least one
element of $R$.

Then the relation spaces of the combinatorial
and canonical quadratic anti-replicators satisfy
\begin{align}
L_2^V(R)
&=
\bigl(\ADup(R)\cap\comp(R)\bigr)
\oplus
\bigoplus_{e_a\notin\supp_{\mathcal B}(R)}
B_{a,2},
\label{eq:comparison-2-revised}\\
L_3^V(R)
&=
\bigl(\ATrip(R)\cap\COMP(R)\bigr)
\oplus
\bigoplus_{e_a\notin\supp_{\mathcal B}(R)}
B_{a,3},
\label{eq:comparison-3-revised}
\end{align}
where $B_{a,2}$ and $B_{a,3}$ are subspaces
of $\mathcal{T}(\gendup)(3)$ and
$\mathcal{T}(\gentrip)(3)$, respectively, with
\[
\dim B_{a,2}=1,\qquad
\dim B_{a,3}=2.
\]
Consequently, the combinatorial
and canonical quadratic anti-replicators coincide if and only if
$R$ has full support with respect to
$\mathcal B$, that is, $|\supp_{\mathcal B}(R)|=3(\dim V)^2$. More precisely, we have
\begin{align}
\dim L_2^V(R)
-\dim\bigl(\ADup(R)\cap\comp(R)\bigr)
&=
3(\dim V)^2-|\supp_{\mathcal B}(R)|,
\label{eq:dimension-2-revised}\\
\dim L_3^V(R)
-\dim\bigl(\ATrip(R)\cap\COMP(R)\bigr)
&=
2\bigl(3(\dim V)^2
-|\supp_{\mathcal B}(R)|\bigr).
\label{eq:dimension-3-revised}
\end{align}
\end{thm}

\begin{proof}
We first determine the additional relations
arising in the canonical quadratic construction.

Fix a monomial $e_a\in\mathcal{T}(V)(3)$.
Up to a permutation of the inputs, we may
assume that its underlying tree is $((xy)z)$.
We use ordered pairs to record the labels
at the root and at the other internal vertex,
respectively. For the anti-bisuccessor, the three marked
expressions are
\begin{equation}\label{eq:three-marked-expressions-revised}
\begin{array}{c|ccc}
 &x&y&z\\ \hline
\ABSu_{\bullet}(e_a)
 &(\prec,\prec)
 &(\prec,\succ)
 &-(\succ,\prec)-(\succ,\succ).
\end{array}
\end{equation}
The four possible split monomials occur
separately, except for the two terms in the
last expression. Thus their difference
belongs to the common kernel of the adjoint
marked maps. This kernel is
\[
B_{a,2}
=
\operatorname{span}
\bigl\{
(\vdash,\dashv)-(\vdash,\vdash)
\bigr\}.
\]
Moreover, since $t^V$ sends each of these two
split monomials to $e_a$, we deduce that $B_{a,2}\subseteq\ker t^V$.

Similarly, for the anti-trisuccessor,
the seven marked expressions are
\begin{equation}\label{eq:seven-marked-expressions-revised}
\begin{array}{c|ccccccc}
J
&\{x\}&\{y\}&\{z\}&\{x,y\}
&\{x,z\}&\{y,z\}&\{x,y,z\}\\ \hline
\ATSu_J(e_a)
&(\prec,\prec)
&(\prec,\succ)
&-\displaystyle\sum_{\epsilon\in\{\prec,\succ,\cdot\}}
(\succ,\epsilon)
&(\prec,\cdot)
&(\cdot,\prec)
&(\cdot,\succ)
&(\cdot,\cdot).
\end{array}
\end{equation}
Hence the corresponding common kernel is $ B_{a,3}=
\operatorname{span}
\left\{
(\vdash,\dashv)-(\vdash,\vdash),\,
(\vdash,\dashv)-(\vdash,\perp)
\right\}. $
We also obtain that $B_{a,3}\subseteq\ker T^V$.
Therefore, we have $\dim B_{a,2}=1$ and $\dim B_{a,3}=2$.

We now compare these two constructions. First, for any $r\in R$ and $x\in\{1,2,3\}$,
choose an element $s_x(r)\in\ADup_x(r)$.
By the adjunction formula
\eqref{eq:adjunction-binary-revised}, we have
\[
q_y^V(s_x(r))=\delta_{xy}r,
\qquad x,y\in\{1,2,3\}.
\]
Moreover, for every $e_a\in\supp_{\mathcal B}(R)$,
there exists $r\in R$ in which $e_a$
has a nonzero coefficient. By changing
the independent choices corresponding to
$e_a$, we obtain
$B_{a,2}\subseteq\ADup(R)$.

Let $h$ be a linear combination of
monomials obtained from the elements of
$\supp_{\mathcal B}(R)$ by replacing each
vertex decoration with either $\dashv$
or $\vdash$. Suppose that
$q_x^V(h)\in R$ for all $x\in\{1,2,3\}$.
Then we have
\[
h-\sum_{x=1}^3s_x(q_x^V(h))
\in\bigoplus_{e_a\in\supp_{\mathcal B}(R)}
B_{a,2}.
\]
Hence $h\in\ADup(R)$.
The converse follows directly from
the definition of the combinatorial
anti-duplicator. Thus, for such $h$,
we obtain that $q_x^V(h)\in R$ for all $x$
if and only if $h\in\ADup(R)$.

Similarly, for every $r\in R$, we can choose
$c(r)\in\comp(R)$ such that $t^V(c(r))=r$.
By the definition of $\comp(R)$, it is
straightforward to verify that
\[
\ker t^V\cap
\operatorname{span}\{e_{a,i}~|~e_a\in\supp_{\mathcal B}(R)\}
\subseteq\comp(R).
\]
Thus, for any linear combination $h$
of the monomials corresponding to
$\supp_{\mathcal B}(R)$, we have
\[
t^V(h)\in R
\quad\Longleftrightarrow\quad
h\in\comp(R).
\]
Therefore, we obtain
\[
L_2^V(R)\big|_{\supp_{\mathcal B}(R)}
=\ADup(R)\cap\comp(R).
\]
Applying the same argument to the
anti-triplicators, we deduce that
\[
L_3^V(R)\big|_{\supp_{\mathcal B}(R)}
=\ATrip(R)\cap\COMP(R).
\]

Now suppose that
$e_a\notin\supp_{\mathcal B}(R)$.
Since no element of $R$ contains $e_a$,
the conditions $q_x^V(h)\in R$ and
$q_J^V(h)\in R$ require the coefficients
of $e_a$ in their images to be zero.
By the preceding computations, the
corresponding relations are exactly
those in $B_{a,2}$ and $B_{a,3}$,
respectively. Moreover,
\[
B_{a,2}\subseteq\ker t^V,
\qquad
B_{a,3}\subseteq\ker T^V.
\]
Since the monomials obtained from
distinct $e_a$ are linearly independent,
we obtain the decomposition formulas
\eqref{eq:comparison-2-revised} and
\eqref{eq:comparison-3-revised}.

Finally, since $\dim B_{a,2}=1$
and $\dim B_{a,3}=2$, the dimension
formulas \eqref{eq:dimension-2-revised}
and \eqref{eq:dimension-3-revised}
follow immediately.
\end{proof}

\subsection{Examples of anti-replicators of binary operads}
We give examples of combinatorial anti-replicators of binary symmetric operads.

\begin{defi}
An {\bf anti-dupli-com algebra} $(A,\diamond)$ is a vector space $A$ equipped with a bilinear operation $\diamond:A \otimes A \to A$ satisfying the following equations for all $x,y,z\in A$:
\begin{eqnarray}
\label{anti-perm-1} (x \diamond y) \diamond z &=& (y \diamond x) \diamond z;\\
\label{anti-perm-2} x \diamond (y \diamond z) &=& y \diamond (x \diamond z);\\
\label{anti-perm-3} (x \diamond y) \diamond z+x \diamond (y \diamond z) &=& (x \diamond z) \diamond y+x \diamond (z \diamond y).
\end{eqnarray}
\end{defi}

\begin{lem}\label{lem:Com-full-support}
The quadratic relation space of the operad $\mathrm{Com}$
has full support with respect to the standard monomial
basis of $\mathcal{T}(V)(3)$.
\end{lem}

\begin{proof}
Let $V=\mathbf{k}\omega$, where $\omega$ is the generating operation of the operad $\com$.
By the convention in \eqref{eq:i-ii-iii}, the space
$\mathcal{T}(V)(3)$ has a basis
\[
e_1=\omega\circ_{\rmi}\omega,\qquad
e_2=\omega\circ_{\rmii}\omega,\qquad
e_3=\omega\circ_{\rmiii}\omega.
\]
The quadratic relation space $R$ of $\mathrm{Com}$
is generated as an $\mathbb{S}_3$-module by
\[
R=\operatorname{span}_{\mathbf{k}}
\{e_1-e_2,\ e_2-e_3\}.
\]
Thus each of the three basis elements $e_1,e_2,e_3$
occurs with a nonzero coefficient in some element
of $R$. Consequently,
\[
\supp_{\mathcal B}(R)=\{e_1,e_2,e_3\},
\]
where $\mathcal B=\{\omega\}$ is the chosen basis of $V$.
Hence $R$ has full support with respect to $\mathcal B$.
\end{proof}

By the above lemma, we have $\widehat{\ADup}(\com)\cong \ADup(\com)$. Thus, the canonical quadratic anti-duplicator and the combinatorial anti-duplicator coincide for the operad $\com$. Henceforth, we shall refer to both constructions simply as the {\bf anti-duplicator of the operad $\com$}.

Denote by $\antiperm$ the operad encoding anti-dupli-com algebras.

\begin{pro}\label{pro:antiperm}
The operad $\antiperm$  is the anti-duplicator of the operad $\com$.
\end{pro}

\begin{proof}
Let $\omega$ be the generating operation of the operad $\com$. Set $\diamond:=\begin{pmatrix}
  \omega\\
  \vdash
  \end{pmatrix}$. Since $\begin{pmatrix}
  \omega\\
  \dashv
  \end{pmatrix}^{(12)}=\begin{pmatrix}
  \omega^{(12)}\\
  \vdash
  \end{pmatrix}=\begin{pmatrix}
  \omega\\
  \vdash
  \end{pmatrix}=\diamond$, we have $\begin{pmatrix}
  \omega\\
  \dashv
  \end{pmatrix}=\diamond^{(12)}.$ The space of relations of $\com$ is generated as an $\mathbb{S}_3$-module by
$$ v_1-v_9= \omega \circ_{\rmi} \omega-\omega \circ_{\rmii} \omega. $$
Then we have
\begin{eqnarray*}
\ADup_{x}(v_1-v_9)&=& \{ z \diamond (y \diamond x)+(y \diamond z) \diamond x,z \diamond (y \diamond x)+(z \diamond y) \diamond x \};\\
\ADup_{y}(v_1-v_9)&=& \{ z \diamond (x \diamond y)-x \diamond (z \diamond y)\};\\
\ADup_{z}(v_1-v_9)&=& \{(x \diamond y) \diamond z-(y \diamond x) \diamond z, (x \diamond y) \diamond z+x \diamond (y \diamond z)\},
\end{eqnarray*}
and
\begin{eqnarray*}
&&\comp(v_1-v_9)= \{(x \diamond y) \diamond z-z \diamond (x \diamond y), (x \diamond y) \diamond z-(y \diamond x) \diamond z, (x \diamond y) \diamond z-z \diamond (y \diamond x),\\
&& (x \diamond y) \diamond z-(y \diamond z) \diamond x,(x \diamond y) \diamond z-x \diamond (y \diamond z),(x \diamond y) \diamond z-(z \diamond y) \diamond x,(x \diamond y) \diamond z-x \diamond (z \diamond y)\}.
\end{eqnarray*}
Computing the displayed quadratic relations yields spanning sets for the relation spaces $\ADup(\com)$ and $\comp(\com)$, respectively, as follows:
\begin{eqnarray*}
\ADup(\com)&=&\langle (x\diamond y)\diamond z- (y\diamond x)\diamond z,~(x\diamond y)\diamond z+x \diamond (y \diamond z) \rangle;\\
\comp(\com)&=&\langle (x\diamond y)\diamond z- (y\diamond x)\diamond z,~ (x\diamond y)\diamond z-x \diamond (y \diamond z),~ (x\diamond y)\diamond z- (x\diamond z)\diamond y \rangle.
\end{eqnarray*}

Now we compute the intersection of subspaces $\ADup(\com)$ and
$\comp(\com)$. It is easy to see that $(x\diamond y)\diamond z-
(y\diamond x)\diamond z$ belongs to the intersection of
$\ADup(\com)$ and $\comp(\com)$, which corresponds to Eq. \eqref{anti-perm-1}. Then we consider the
following equation
\begin{eqnarray*}
&&a_1 ((x\diamond y)\diamond z+x \diamond (y \diamond z))+a_2 ((x\diamond z)\diamond y+x \diamond (z \diamond y))+ a_3 (y \diamond (x \diamond z)+(y \diamond x)\diamond z)\\
&&+a_4 (y \diamond (z \diamond x)+(y \diamond z)\diamond x)+a_5 ((z \diamond x)\diamond y+ z \diamond (x \diamond y))+a_6 ((z \diamond y)\diamond x+z \diamond (y \diamond x))\\
&=& b_1 ((x \diamond y)\diamond z-(y \diamond x)\diamond z)+b_2  ((x \diamond z)\diamond y-(z \diamond x)\diamond y)+ b_3 ((y \diamond z)\diamond x-(z \diamond y)\diamond x)\\
&&+ b_4 ((x \diamond y)\diamond z-x \diamond (y \diamond z))+b_5 ((x \diamond z)\diamond y-x \diamond (z \diamond y))+b_6 ((y \diamond x)\diamond z-y \diamond (x \diamond z))\\
&&+b_7 ((y \diamond z)\diamond x-y \diamond (z \diamond x))+b_8 ((z \diamond x)\diamond y-z \diamond (x \diamond y))+ b_9 ((z \diamond y)\diamond x-z \diamond (y \diamond x))\\
&&+b_{10} ((x \diamond y)\diamond z-(x \diamond z) \diamond y)+b_{11} ((y \diamond z)\diamond x-(y \diamond x)\diamond z)+b_{12} ((z \diamond x)\diamond y-(z \diamond y) \diamond x).
\end{eqnarray*}
Since $\huaT(V)(3)$ is generated by the twelve elements $v_i~(1\leq i \leq 12)$ in the table \eqref{eq:table}, we obtain the following system of linear equations in the variables $a_m, b_n$, where $1 \leq m \leq 6$ and $1 \leq n \leq 12$.
$$ \begin{array}{c@{\quad}c@{\quad}c}
\left\{
\begin{aligned}
a_1 &= -b_4, \\
a_2 &= -b_5, \\
a_3 &= -b_6, \\
a_4 &= -b_7, \\
a_5 &= -b_8, \\
a_6 &= -b_9,
\end{aligned}
\right.
&
\text{and}
&
\left\{
\begin{aligned}
2a_1-b_1-b_{10} &= 0, \\
2a_5+b_2-b_{12} &= 0, \\
2a_4-b_3-b_{11} &= 0, \\
2a_2-b_2+b_{10} &= 0, \\
2a_6+b_3+b_{12} &= 0, \\
2a_3+b_1+b_{11} &= 0.
\end{aligned}
\right.
\end{array} $$
The coefficient matrix of this system has rank $12$. Hence its solution space is $6$-dimensional. Thus the general solution of the above system of linear equations can be expressed as
$$
(a_1,\dots,a_6;b_1,\dots,b_{12})^{\mathrm{T}}
=
t_1\eta_1+t_2\eta_2+t_3\eta_3+t_4\eta_4+t_5\eta_5+t_6\eta_6,
$$
where $t_1,\dots,t_6\in \mathbf{k}$ are arbitrary constants.

A basis of the solution space is given by $(\eta_1,\dots,\eta_6)$, where
\begin{eqnarray*}
\begin{aligned}
\eta_1&=(1,0,-1,0,0,0;\ 2,0,0,-1,0,1,0,0,0,0,0,0);\\
\eta_2&=(0,1,0,0,-1,0;\ 0,2,0,0,-1,0,0,1,0,0,0,0);\\
\eta_3&=(0,0,0,1,0,-1;\ 0,0,2,0,0,0,-1,0,1,0,0,0);\\
\eta_4&=(1,-1,0,0,0,0;\ 0,0,0,-1,1,0,0,0,0,2,0,0);\\
\eta_5&=(-1,0,0,1,0,0;\ -2,0,0,1,0,0,-1,0,0,0,2,0);\\
\eta_6&=(0,1,0,-1,0,0;\ 0,2,-2,0,-1,0,1,0,0,0,0,2).
\end{aligned}
\end{eqnarray*}
Here the coordinate order is
$ (a_1,a_2,a_3,a_4,a_5,a_6;b_1,b_2,b_3,b_4,b_5,b_6,b_7,b_8,b_9,b_{10},b_{11},b_{12}). $
The six vectors $\eta_1,\dots,\eta_6$ are linearly independent and hence form a basis of the solution space of the above system.

The solution $\eta_1$ yields the relation $(x \diamond y)\diamond z-(y \diamond x)\diamond z+x \diamond (y \diamond z)-y \diamond (x \diamond z).$ Since $(x\diamond y)\diamond z- (y\diamond x)\diamond z$ belongs to the intersection of $\ADup(\com)$ and $\comp(\com)$, we deduce that the relation
$x\diamond (y \diamond z)-y \diamond (x \diamond z)$ belongs to the intersection of $\ADup(\com)$ and $\comp(\com)$, which corresponds to Eq. \eqref{anti-perm-2}. The solutions $\eta_2$ and $\eta_3$ yield analogous relations up to permutation of the variables.

Considering the solution $\eta_4$, we have the relation $(x \diamond y)\diamond z-(x \diamond z)\diamond y+x\diamond (y \diamond z)-x\diamond (z \diamond y)$, which corresponds to Eq. \eqref{anti-perm-3}. The solutions $\eta_5$ and $\eta_6$ yield analogous relations up to permutation of the variables. Therefore, the operad $\antiperm$ is the anti-duplicator of the operad $\com$.
\end{proof}

\begin{defi}
An {\bf anti-dupli-Lie algebra} $(A,\circ)$ is a vector space $A$ equipped with a bilinear operation $\circ:A \otimes A \to A$ satisfying the following equations
\begin{eqnarray}
\label{anti-leibniz-1} (x \circ y+y \circ x)\circ z&=&0;\\
\label{anti-leibniz-2} \sum_{\text{cyclic}} (x \circ (z \circ y)+(x \circ z)\circ y-z \circ (x \circ y))&=&0, \quad \forall x,y,z \in A,
\end{eqnarray}
where the sum is taken over cyclic permutations of $x,y,z \in A$.
\end{defi}

\begin{lem}\label{lem:Lie-full-support}
The quadratic relation space of the operad $\mathrm{Lie}$
has full support with respect to the standard monomial
basis of $\mathcal{T}(V)(3)$.
\end{lem}

\begin{proof}
Let $V=\mathbf{k}\mu$, where $\mu$ is the generating operation of the operad $\lie$. By the convention in  \eqref{eq:i-ii-iii}, the space $\mathcal{T}(V)(3)$ has a basis
\[
e_1=\mu\circ_{\rmi}\mu,\qquad
e_2=\mu\circ_{\rmii}\mu,\qquad
e_3=\mu\circ_{\rmiii}\mu.
\]
The quadratic relation space $R$ of $\mathrm{Lie}$ is generated by $e_1+e_2+e_3=0$. Hence
\[
R=\operatorname{span}_{\mathbf{k}}
\{e_1+e_2+e_3\}.
\]
Since all three coefficients are nonzero, we have
\[
\supp_{\mathcal B}(R)=\{e_1,e_2,e_3\},
\]
where $\mathcal B=\{\mu\}$ is the chosen basis of $V$.
Therefore, $R$ has full support with respect to $\mathcal B$.
\end{proof}

By the above lemma, we have $\widehat{\ADup}(\lie)\cong \ADup(\lie)$. Thus, the canonical quadratic anti-duplicator and the combinatorial anti-duplicator coincide for the operad $\lie$. Henceforth, we shall refer to both constructions simply as the {\bf anti-duplicator of the operad $\lie$}.

Denote by $\antileibniz$ the operad encoding anti-dupli-Lie algebras.

\begin{pro}\label{pro:antileibniz}
The operad $\antileibniz$  is the anti-duplicator of the operad $\lie$.
\end{pro}

\begin{proof}
Let $\mu$ be the generating operation of the operad $\lie$. Set $\circ:=\begin{pmatrix}
  \mu\\
  \vdash
  \end{pmatrix}$. Since $\begin{pmatrix}
  \mu\\
  \dashv
  \end{pmatrix}^{(12)}=\begin{pmatrix}
  \mu^{(12)}\\
  \vdash
  \end{pmatrix}=-\begin{pmatrix}
  \mu\\
  \vdash
  \end{pmatrix}=-~\circ$, we have $\begin{pmatrix}
  \mu\\
  \dashv
  \end{pmatrix}=-\circ^{(12)}.$ The space of the relations of $\lie$ is generated as an $\mathbb{S}_3$-module by
$$ v_1+v_5+v_9= \mu \circ_{\rmi} \mu+\mu \circ_{\rmii} \mu+\mu \circ_{\rmiii} \mu. $$
Then we have
\begin{eqnarray*}
\ADup_{x}(v_1+v_5+v_9)&=& \{z \circ (y \circ x)-y \circ (z \circ x)-(y \circ z)\circ x,z \circ (y \circ x)-y \circ (z \circ x)+(z \circ y) \circ x\};\\
\ADup_{y}(v_1+v_5+v_9)&=& \{x \circ (z \circ y)-z \circ (x \circ y)-(z \circ x)\circ y,x \circ (z \circ y)-z \circ (x \circ y)+(x \circ z) \circ y\};\\
\ADup_{z}(v_1+v_5+v_9)&=& \{y \circ (x \circ z)-x \circ (y \circ z)-(x \circ y)\circ z,y \circ (x \circ z)-x \circ (y \circ z)+(y \circ x) \circ z\}.
\end{eqnarray*}
It is straightforward to check that the relation space $\ADup(\lie)$ is spanned by
\begin{eqnarray}\label{eq:ADup-Lie}
\ADup(\lie)=\langle x \circ (y \circ z)+(x \circ y)\circ z-y \circ
(x \circ z)\rangle.
\end{eqnarray}
A direct calculation also shows that the relation space $\comp(\lie)$ is spanned by
\begin{eqnarray*}
\comp(\lie)=\langle (x \circ y) \circ z+(z \circ x) \circ y+ (y \circ z) \circ x,
 ~(x \circ y+y \circ x)\circ z,
 ~(x \circ y) \circ z+z \circ (x \circ y) \rangle.
\end{eqnarray*}

Now we compute the intersection of subspaces $\ADup(\lie)$ and
$\comp(\lie)$. Since
$$(y \circ z+z \circ y)\circ x=(x\circ y)\circ z+(z\circ x)\circ y+(y\circ z)\circ x+(y\circ x)\circ z+(z \circ y)\circ x+(x \circ z)\circ y-(x \circ y+y \circ x)\circ z-(z \circ x+x \circ z)\circ y, $$
we consider the following equation:
\begin{eqnarray*}
&&a_1 ((x\circ y)\circ z+(z\circ x)\circ y+(y\circ z)\circ x)+a_2 ((y\circ x)\circ z+(z \circ y)\circ x+(x \circ z)\circ y)\\
&&+ a_3 ((x \circ y)\circ z+ z\circ (x \circ y))+a_4 ((y \circ x)\circ z+ z\circ (y \circ x))+a_5 ((z \circ y)\circ x+ x\circ (z \circ y))\\
&&+a_6 ((x \circ z)\circ y+ y \circ (x \circ z))+a_7 ((y \circ z)\circ x+ x \circ (y \circ z))+a_8 ((z \circ x)\circ y+ y \circ (z \circ x))\\
&&+a_9 ((x \circ y+y \circ x)\circ z)+a_{10} ((z \circ x+x \circ z)\circ y)\\
&=& b_1 (x \circ (y \circ z)+(x \circ y)\circ z-y \circ (x\circ z))+b_2  (x \circ (z \circ y)+(x \circ z)\circ y-z \circ (x\circ y))\\
&&+b_3 (y \circ (x \circ z)+(y \circ x)\circ z-x \circ (y\circ z))+b_4  (z \circ (y \circ x)+(z \circ y)\circ x-y \circ (z\circ x))\\
&&+b_5 (y \circ (z \circ x)+(y \circ z)\circ x-z \circ (y\circ x))+b_6  (z \circ (x \circ y)+(z \circ x)\circ y-x \circ (z\circ y)).
\end{eqnarray*}
Since $\huaT(V)(3)$ is generated by the twelve elements $v_i~(1\leq i \leq 12)$ in the table \eqref{eq:table}, we obtain the following system of linear equations in the variables $a_m$ and $b_n$, where $1 \leq m \leq 10$ and $1 \leq n \leq 6$,
$$ \begin{array}{c@{\quad}c@{\quad}c}
\left\{
\begin{aligned}
a_1+a_3+a_9 &= b_1, \\
a_1+a_8+a_{10} &= b_6, \\
a_1+a_7 &= b_5, \\
a_2+a_6+a_{10} &= b_2, \\
a_2+a_5 &= b_4, \\
a_2+a_4+a_9 &= b_3,
\end{aligned}
\right.
&
\text{and}
&
\left\{
\begin{aligned}
a_7-b_1+b_3 &= 0, \\
a_3+b_2-b_6 &= 0, \\
a_8+b_4-b_5 &= 0, \\
a_5-b_2+b_6 &= 0, \\
a_4-b_4+b_5 &= 0, \\
a_6+b_1-b_3 &= 0.
\end{aligned}
\right.
\end{array} $$
The coefficient matrix of this system has rank $12$. Hence its solution space is $4$-dimensional. Thus the general solution of the above system of linear equations can be expressed as
$$
(a_1,\dots,a_{10};b_1,\dots,b_6)^{\mathrm{T}}
=
s_1\delta_1+s_2\delta_2+s_3\delta_3+s_4\delta_4,
$$
where $s_1,\dots,s_4\in \mathbf{k}$ are arbitrary constants.

A basis of the solution space is given by $(\delta_1,\delta_2,\delta_3,\delta_4)$, where
\begin{eqnarray*}
\begin{aligned}
\delta_1&=(1,1,0,0,0,0,0,0,0,0;\ 1,1,1,1,1,1);\\
\delta_2&=(1,0,-1,1,1,1,-1,-1,0,0;\ 0,1,1,1,0,0);\\
\delta_3&=(0,0,0,0,0,0,0,0,1,0;\ 1,0,1,0,0,0);\\
\delta_4&=(0,0,0,0,0,0,0,0,0,1;\ 0,1,0,0,0,1).
\end{aligned}
\end{eqnarray*}
Here the coordinate order is
$ (a_1,a_2,a_3,a_4,a_5,a_6,a_7,a_8,a_9,a_{10};b_1,b_2,b_3,b_4,b_5,b_6). $
The four vectors $\delta_1,\dots,\delta_4$ are linearly independent and hence form a basis of the solution space of the above system.

The solutions $\delta_3$, $\delta_1-\delta_3-\delta_4$ and $\delta_4$
yield the relations
\[
(x\circ y+y\circ x)\circ z,\quad
(y\circ z+z\circ y)\circ x,\quad
(z\circ x+x\circ z)\circ y,
\]
respectively, which correspond to
Eq.~\eqref{anti-leibniz-1}.
On the other hand, the solution $\delta_2$
gives rise to the relation
\[
\sum_{\mathrm{cyclic}}
\bigl(x\circ(z\circ y)+(x\circ z)\circ y
-z\circ(x\circ y)\bigr),
\]
which corresponds to Eq.~\eqref{anti-leibniz-2}.
Consequently, the operad $\antileibniz$ is the anti-duplicator of the operad $\lie$.
\end{proof}

Let us compare the canonical quadratic anti-duplicator with the combinatorial anti-duplicator of the operad $\mathrm{2Ass}$.

\begin{ex}\label{ex:2Ass-canonical}
Let $\mathrm{2Ass}=\mathcal{T}(V)/(R)$ be the operad
of $2$-associative algebras~\cite{Lo08,LR2}, generated by two
associative operations $\cdot_1$ and $\cdot_2$
without further relations.
Here $V$ has the compatible basis
\[
\mathcal B=
\{\mu_1,\mu_1^{(12)},\mu_2,\mu_2^{(12)}\}.
\]
We use the induced monomial basis of $\huaT(V)(3)$.
The relation space $R$ is generated as an $\mathbb{S}_3$-module
by the associativity relations
\[
(x\cdot_i y)\cdot_i z-x\cdot_i(y\cdot_i z),
\qquad i=1,2.
\]
Write $\dashv_i$ and $\vdash_i$ for the two
operations corresponding to $\mu_i$ in $\gendup$. For each $i=1,2$, the combinatorial
anti-duplicator gives the following relations:
\begin{align}
x\dashv_i(y\dashv_i z)
&=x\dashv_i(y\vdash_i z),
\label{eq:2Ass-pure-1}\\
(x\dashv_i y)\vdash_i z
&=(x\vdash_i y)\vdash_i z,
\label{eq:2Ass-pure-2}\\
(x\vdash_i y)\dashv_i z
&=x\vdash_i(y\dashv_i z),
\label{eq:2Ass-pure-3}\\
(x\dashv_i y)\dashv_i z
+x\dashv_i(y\dashv_i z)
&=(x\dashv_i y)\vdash_i z
+x\vdash_i(y\vdash_i z).
\label{eq:2Ass-pure-4}
\end{align}
There are $24$ independent relations for each $i$, hence
$48$ pure quadratic relations in total.

On the other hand,
$\dim\mathcal{T}(V)(3)=48$, whereas
$|\supp_{\mathcal B}(R)|=24$.
By Theorem~\ref{lem:comparison-two-rep},
the remaining $24$ mixed monomials contribute
one additional relation each through $B_{a,2}$.
These relations are generated as an
$\mathbb{S}_3$-module by
\begin{align}
(x\dashv_j y)\vdash_i z
&=(x\vdash_j y)\vdash_i z,
\label{eq:2Ass-mixed-1}\\
x\dashv_i(y\dashv_j z)
&=x\dashv_i(y\vdash_j z),
\label{eq:2Ass-mixed-2}
\end{align}
where $i,j\in\{1,2\}$ and $i\neq j$.

Consequently, the relations
\eqref{eq:2Ass-pure-1}--\eqref{eq:2Ass-mixed-2}
generate the entire quadratic relation space of
$\widehat{\ADup}(\mathrm{2Ass})$. Moreover, we deduce that $\dim L_2^V(R)=48+24=72.$
\end{ex}

\begin{rmk}\label{rmk:2Ass-comparison}
The difference between
$\ADup(\mathrm{2Ass})$ and
$\widehat{\ADup}(\mathrm{2Ass})$
lies precisely in the mixed quadratic relations.
The combinatorial construction gives only
the $48$ pure relations, whereas the canonical
quadratic construction introduces $24$
additional relations arising from the spaces
$B_{a,2}$ associated with monomials outside
$\supp_{\mathcal B}(R)$.
\end{rmk}

\section{Koszul duality between anti-successors and canonical quadratic anti-replicators}
\label{sec:dual}

We establish the Koszul duality between
anti-successors and canonical quadratic anti-replicators of binary quadratic operads. We then show that the canonical
quadratic anti-duplicator of a binary quadratic operad $\opd$ is isomorphic to the Manin white product of $\antiperm$ and $\opd$.

Let $\opd=\mathcal{T}(\gensp)/(\relsp)$ be a binary quadratic
operad with a finite-dimensional
$\mathbb{S}_2$-module $\gensp$ and
$\relsp\subseteq\mathcal{T}(\gensp)(3)$ be an
$\mathbb{S}_3$-submodule. Recall that
\[
\gensp^\vee=\gensp^*\otimes\sgn_2,
\qquad
\opd^!=\mathcal{T}(\gensp^\vee)/(\relsp^\perp).
\]
By the definitions of the generating spaces, we have
natural isomorphisms of $\mathbb{S}_2$-modules
\[
(\genbi)^\vee\cong\overline{\gensp^\vee},
\qquad
(\gentri)^\vee\cong\overrightarrow{\gensp^\vee}.
\]
More precisely, for $\omega\in\gensp$ and
$f\in\gensp^\vee$, these identifications are given by
\[
\begin{gathered}
\left\langle
\begin{pmatrix}f\\ \dashv\end{pmatrix},
\begin{pmatrix}\omega\\ \prec\end{pmatrix}
\right\rangle
=
\langle f,\omega\rangle,
\qquad
\left\langle
\begin{pmatrix}f\\ \vdash\end{pmatrix},
\begin{pmatrix}\omega\\ \succ\end{pmatrix}
\right\rangle
=
\langle f,\omega\rangle,
\qquad
\left\langle
\begin{pmatrix}f\\ \perp\end{pmatrix},
\begin{pmatrix}\omega\\ \cdot\end{pmatrix}
\right\rangle
=
\langle f,\omega\rangle.
\end{gathered}
\]
where the pairings between different labels vanish.
Here $\langle f,\omega\rangle$ denotes the natural
pairing between $\gensp^\vee$ and $\gensp$.

The following relations provide the key identities for
the Koszul duality.

\begin{lem}\label{lem:koszul-duality}\label{lem:annihilator-revised}
With the above notations, we have
\begin{align}
\bigl(\ABSu(\relsp)+\st(\relsp)\bigr)^\perp
&=L_2^{\gensp^\vee}(\relsp^\perp),
\label{eq:orthogonality-2-revised}\\
\bigl(\ATSu(\relsp)+\ST(\relsp)\bigr)^\perp
&=L_3^{\gensp^\vee}(\relsp^\perp).
\label{eq:orthogonality-3-revised}
\end{align}
Here $\ABSu(\relsp)$ and $\ATSu(\relsp)$
denote the sums of the images of
$\ABSu_x$ for $x=1,2,3$ and $\ATSu_J$
for all nonempty subsets $J\subseteq\{1,2,3\}$,
respectively.
\end{lem}

\begin{proof}
Let $h\in\mathcal{T}(\overline{\gensp^\vee})(3)$.
By \eqref{eq:adjunction-binary-revised},
for every $r\in\relsp$, we have
\[
\langle h,\ABSu_x(r)\rangle
=\langle q_x^{\gensp^\vee}(h),r\rangle,
\qquad
\langle h,\st(r)\rangle
=\langle t^{\gensp^\vee}(h),r\rangle.
\]
Hence $h$ is orthogonal to
$\ABSu(\relsp)+\st(\relsp)$ if and only if
\[
q_x^{\gensp^\vee}(h)\in\relsp^\perp
\quad(x=1,2,3),\qquad
t^{\gensp^\vee}(h)\in\relsp^\perp.
\]
By Definition~\ref{def-can-rep}, these are
precisely the conditions defining
$L_2^{\gensp^\vee}(\relsp^\perp)$.
Hence \eqref{eq:orthogonality-2-revised} holds.

Similarly, for
$h\in\mathcal{T}(\overrightarrow{\gensp^\vee})(3)$,
replacing
$\ABSu_x,\st,q_x^{\gensp^\vee},t^{\gensp^\vee}$
by
$\ATSu_J,\ST,q_J^{\gensp^\vee},T^{\gensp^\vee}$,
respectively, we deduce that
\eqref{eq:orthogonality-3-revised} holds.
\end{proof}

\begin{thm}\label{thm:koszul-dual}
Let $\opd=\mathcal{T}(\gensp)/(\relsp)$ be a binary quadratic operad with a
finite-dimensional generating $\mathbb{S}_2$-module $V$.
Then there is a natural isomorphism of operads
\[
\widehat{\ADup}(\opd^!)\cong\ABSu(\opd)^!.
\]
\end{thm}

\begin{proof}
Since $\opd^!=\mathcal{T}(\gensp^\vee)/(\relsp^\perp)$ and the relation space of $\ABSu(\opd)$ is $\ABSu(\relsp)+\st(\relsp)$, by \eqref{eq:orthogonality-2-revised}, we have
\[
\begin{aligned}
\ABSu(\opd)^!
&=\mathcal{T}(\genbi^\vee)
\big/
\bigl((\ABSu(\relsp)+\st(\relsp))^\perp\bigr) \cong
\mathcal{T}(\overline{\gensp^\vee})
\big/
\bigl((\ABSu(\relsp)+\st(\relsp))^\perp\bigr)\\
&=
\mathcal{T}(\overline{\gensp^\vee})
\big/
\bigl(L_2^{\gensp^\vee}(\relsp^\perp)\bigr) =
\widehat{\ADup}(\opd^!).
\end{aligned}
\]
This completes the proof.
\end{proof}

\begin{thm}\label{thm:koszul-dual-trip}
Let $\opd=\mathcal{T}(\gensp)/(\relsp)$ be a binary quadratic operad with a
finite-dimensional generating $\mathbb{S}_2$-module $V$.
Then there is a natural isomorphism of operads
\[
\widehat{\ATrip}(\opd^!)\cong\ATSu(\opd)^!.
\]
\end{thm}

\begin{proof}
Since the relation space of $\ATSu(\opd)$ is $\ATSu(\relsp)+\ST(\relsp)$, by Eq.~\eqref{eq:orthogonality-3-revised},
we obtain
\[
\begin{aligned}
\ATSu(\opd)^!
&= \mathcal{T}(\gentri^\vee)
\big/\bigl((\ATSu(\relsp)+\ST(\relsp))^\perp\bigr)
 \cong
\mathcal{T}(\overrightarrow{\gensp^\vee})
\big/
\bigl((\ATSu(\relsp)+\ST(\relsp))^\perp\bigr)\\
&=
\mathcal{T}(\overrightarrow{\gensp^\vee})
\big/
\bigl(L_3^{\gensp^\vee}(\relsp^\perp)\bigr) =
\widehat{\ATrip}(\opd^!).
\end{aligned}
\]
This completes the proof.
\end{proof}

As consequences, we obtain the following Koszuldual identifications.

\begin{cor}\label{cor:perm-prelie}
\begin{enumerate}
\item \label{i:pp1} The operad $\antiperm$ is the Koszul dual of the operad $\antiprelie$.
\item
% [11/04/2026] Cross-reference anchor only; corollary itself unchanged.
\label{i:pp2}
The operad $\antileibniz$ is the Koszul dual of the operad $\antizinbiel$.
\end{enumerate}
\end{cor}

\begin{proof}
\eqref{i:pp1}
By Theorem \ref{thm:koszul-dual} and Propositions~\ref{pro:antiperm} and
\ref{pro:anti-pre-lie}, with $\opd:=\lie$, we have
$$ \antiperm=\widehat{\ADup}(\com)=\widehat{\ADup}(\lie^!)=\ABSu(\lie)^{!}=\antiprelie^{!}, $$
which implies that the operad $\antiperm$ is the Koszul dual of the operad $\antiprelie$.

\noindent
\eqref{i:pp2}
This follows from the same argument.
\end{proof}

We conclude by applying the preceding duality theorem to describe the Manin white product of $\antiperm$ with $\opd$ in terms of canonical anti-duplicators. For background on Manin white products, see \cite{GK2014,Val2}.

Let $\opd=\mathcal{T}(\gensp)/ (\relsp)$ and $\oqd=\mathcal{T}(W)/ (S)$ be two binary quadratic operads with finite-dimensional generating spaces. For the operad morphism ~\cite{GK95}
$$\Phi: \mathcal{T}(\gensp\otimes W) \to \mathcal{T}(\gensp)\otimes \mathcal{T}(W),$$ consider the composition
$$ \mathcal{T}(\gensp\otimes W) \overset{\Phi}{\longrightarrow}  \mathcal{T}(\gensp)\otimes \mathcal{T}(W)  \overset{\pi_{\opd}\otimes \pi_{\oqd}}{\longrightarrow}
\opd \otimes \oqd, $$
where $\pi_{\opd}:\mathcal{T}(\gensp) \to \opd$ and  $\pi_{\oqd}:\mathcal{T}(W) \to \oqd$ are natural projections.
Let
$$
\Phi_3^{-1}\Bigl(
     R\otimes\mathcal T(W)(3)
     +\mathcal T(V)(3)\otimes S\Bigr):=
\Phi^{-1}\Bigl(
     R\otimes\mathcal T(W)(3)
     +\mathcal T(V)(3)\otimes S\Bigr) \cap
\mathcal T(V\otimes W)(3)
$$
be the quadratic component of the kernel.

\begin{defi}
Let $\opd=\mathcal{T}(\gensp)/ (\relsp)$ and $\oqd=\mathcal{T}(W)/(S)$ be two binary quadratic operads with finite-dimensional generating spaces. The {\bf (quadratic) Manin white product} of $\opd$ and $\oqd$ is defined by
$$ \opd \bigcirc  \oqd :=\mathcal{T}(\gensp\otimes W)/ (\Phi^{-1}_3(R\otimes \mathcal{T}(W)+\mathcal{T}(\gensp)\otimes S)). $$
\end{defi}

In general, it is difficult to compute the Manin white product
when the operads are given in terms of generators and relations.
The following theorem provides an efficient way to
compute the white Manin product of a binary quadratic operad with
the operad $\antiperm$ by relating them to the anti-duplicator.

\begin{thm}\label{thm:ADup-antiperm}
Let $\opd=\mathcal{T}(\gensp)/ (\relsp)$ be a binary quadratic operad with a finite-dimensional generating $\mathbb{S}_2$-module. Then we have the isomorphism of operads
$$  \widehat{\ADup}(\opd)\cong \antiperm  \bigcirc \opd.$$
\end{thm}

\begin{proof}
By Theorems \ref{thm:Manin-product}, \ref{thm:koszul-dual} and
Corollary \ref{cor:perm-prelie}, we have
$$  \widehat{\ADup}(\opd)=\ABSu(\opd^!)^! \cong (\antiprelie \bullet \opd^!)^! \cong \antiprelie^! \bigcirc (\opd^!)^! \cong \antiperm \bigcirc \opd.$$
Here we use the fact that $(\opd \bullet \oqd)^! \cong \opd^! \bigcirc \oqd^!$ \cite[Theorem 8.8.8]{LV}. This completes the proof.
\end{proof}

Applying Theorem \ref{thm:ADup-antiperm} with $\opd=\lie$ and using Proposition \ref{pro:antileibniz},
we immediately get the following consequence.

\begin{cor}
The Manin white product of the operad $\antiperm$ and the operad $\lie$ is isomorphic to the operad $\antileibniz$
$$\antiperm  \bigcirc \lie \cong \antileibniz.$$
\end{cor}

\noindent{\bf Acknowledgments.} We express our sincerest gratitude to Vladimir Dotsenko for his valuable comments on computing the intersection of the relation spaces of binary operads. This work is supported by NSFC
 (12271265, 12471060, W2412041, 12671037, 12601051),
the Central Government-Guided Fund Project for Local Science and Technology Development
in Tianjin (25ZYZYJC00010), the Fundamental Research Funds for the Central Universities,
and the Nankai Zhide Foundation.

\noindent
{\bf Declaration of interests. } The authors have no conflict of interest to declare.

\noindent
{\bf Data availability. }
Data sharing is not applicable as no data were created or analyzed.

\end{document}